\documentclass[11pt]{article}

\usepackage{amsmath,amssymb,amscd,amsthm,amsfonts}
\usepackage{graphicx,subfigure}
\usepackage{hyperref}
\usepackage{dsfont}
\usepackage{color}
\usepackage{float}

\newtheorem{theorem}{Theorem}[section]

\newtheorem{corollary}[theorem]{Corollary}
\newtheorem{conjecture}[theorem]{Conjecture}

\newtheorem{problem}[theorem]{Problem}

\newcommand{\N}{\mathbb{N}}
\newcommand{\rr}{\mathds{R}}
\newcommand{\R}{\mathds{R}}
\newcommand{\C}{\mathcal{C}}
\newcommand{\Z}{\mathbb{Z}}

\newcommand{\p}{{\mathcal P}}

\newcommand{\tri}{\Delta}
\newcommand{\al}{\alpha}
\newcommand{\si}{\sigma}
\def\rr{\mathds{R}}
\DeclareMathOperator{\conv}{conv}
\DeclareMathOperator{\skel}{skel}
\DeclareMathOperator{\relint}{relint}
\DeclareMathOperator{\dist}{dist}

\DeclareMathOperator{\aff}{aff}

\title{Tverberg's theorem is 50 years old: a survey}

\author{Imre B\'ar\'any \and Pablo Sober\'on}

\begin{document}

\maketitle

\begin{abstract}
This survey presents an overview of the advances around Tverberg's theorem, focusing on the last two decades.  We discuss the topological, linear-algebraic, and combinatorial aspects of Tverberg's theorem and its applications.  The survey contains several open problems and conjectures.
\end{abstract}

\section{Introduction}\label{section-introduction}

Tverberg's theorem has been a cornerstone of combinatorial convexity for over fifty years. Its impact and influence is only comparable to that of the famous and classic theorems of Carath\'eodory and Helly. This gem lies at the crossroads of combinatorics, topology, and linear algebra, and continues to yield challenging and interesting open problems.  Its states the following.

\begin{theorem}[Helge Tverberg 1966 \cite{Tverberg:1966tb}]
Given $(r-1)(d+1)+1$ points in $\rr^d$, there is a partition of them into $r$ parts whose convex hulls intersect.
\end{theorem}

\begin{figure}
\centerline{\includegraphics[scale=1]{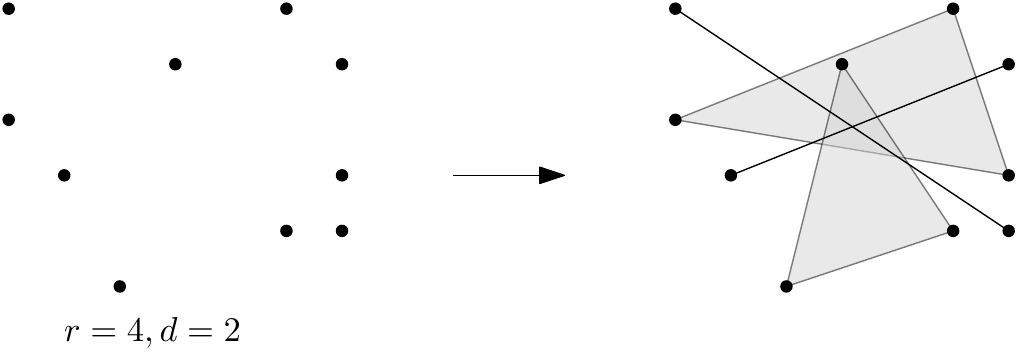}}
\caption{An example of a Tverberg partition.  The partition is not unique.}
\end{figure}

More formally, given $X\subset \rr^d$ of $(r-1)(d+1)+1$ points, there is a partition $X=X_1\cup \dots \cup X_r$ such that $\bigcap_{j=1}^r \conv X_j \ne \emptyset$. Such a partition is called a
 \textit{Tverberg partition}. The number of points in this result is optimal, as a dimension-counting argument shows. In fact, if $X$ is in general enough position and in the partition $X=X_1\cup \ldots \cup X_r$ we have $1\le |X_j|\le d+1$ for every $j$, then $\bigcap_{j=1}^r \aff X_j$ is a single point if $|X|= (r-1)(d+1)+1$, and is empty if $|X|\le (r-1)(d+1)$.


The last decade has seen an impressive sequence of results around Tverberg's theorem.  The purpose of this survey is to give a broad overview of  the current state of the field and point out key open problems.  Other surveys covering different aspects of Tverberg's theorem can be found in \cite{Eckhoff:1979bi, Eck93survey, Matousek:2002td, BBZ17survey, de2017discrete, BZ17}.

   The paper is organized as follows.  In sections \ref{section-topological} and \ref{section-colored} we describe the topological and colorful versions of Tverberg's theorem, which have received the most attention in recent years.  In sections \ref{section-intersection} and \ref{section-universal} we discuss a large number of variations and conjectures around Tverberg's theorem.  In Section \ref{section-applications} we describe some applications of Tverberg's theorem.  Finally, in Section \ref{section-spaces} we present Tverberg-type results where the settings have changed dramatically, such as Tverberg for convexity spaces or quantitative versions.  In that last section, we focus mostly on results which are related to geometry.

\subsection{Interlude: a short history of Tverberg's theorem}
An early predecessor of Tverberg's theorem is Radon's lemma from 1921 \cite{Radon:1921vh, Eckhoff:1979bi}. Radon used it in his proof of Helly's theorem. It says that \textit{any set $X$ of $d+2$ points in $\rr^d$ can be split into two sets whose convex hulls intersect}. So it is the case $r=2$ of Tverberg's theorem. Its proof is simple: the $d+2$ vectors in $X$ have a nontrivial affine dependence $\sum_{x \in X}\al(x)x=0$ and  $\sum_{x \in X}\al(x)=0$. The sets $X_1=\{x \in X: \al(x)\ge 0\}$ and $X_2=\{x \in X: \al(x) < 0\}$ form a partition of $X$ and their convex hulls intersect, as one can easily check.

Another result linked to this theorem is Rado's centerpoint theorem.  This states that \textit{for any set $X$ of $n$ points in $\rr^d$, there is a point $p$ such that any closed half-space that contains $p$ also contains at least $\left\lceil \frac{n}{d+1}\right\rceil$ points of $X$}. The standard proof of this result uses Helly's theorem. Tverberg's theorem implies it in few lines: setting $r=\left\lceil \frac{n}{d+1}\right\rceil$, there is a partition of $X$ into $r$ parts $X_1,\ldots,X_r$ and a point $p\in \rr^d$ such that $p \in \bigcap_{j=1}^r \conv X_j$. Then $p$ is a centerpoint of $X$: every closed halfspace containing $p$ contains at least one point from each $X_j$.

In a paper entitled ``On $3N$ points in a plane'' Birch~\cite{Birch:1959} proves that any $3N$ points in the plane determine $N$ triangles that have a point in common. His motivation was the (planar) centerpoint theorem. Actually, he proves more, namely the case $d=2$ of Tverberg's theorem and states the general case as a conjecture.

Tverberg's original motivation was also the centerpoint theorem and he learned about Birch's result and conjecture only later. He proved it first for $d=3$ in 1963, and in full generality in 1964.  Here is, in his own words, how he found the proof:  ``I recall that the weather was bitterly cold in Manchester. I awoke very early one morning shivering, as the electric heater in the hotel room had gone off, and I did not have an extra shilling to feed the meter. So, instead of falling back to sleep, I reviewed the problem once more, and then the solution dawned on me!'' \cite{tve:recollections}.


\subsection{Proof methods}

By now there are several proofs of Tverberg's theorem, two by Tverberg himself \cite{Tverberg:1966tb, Tve81}, one by Tverberg and Vre\'cica \cite{TV93}, by Roudneff \cite{Roudneff:2001cl}, by Sarkaria \cite{Sarkaria:1992vt}, and by Zvagelskii \cite{Zva08}. We explain here two of them. The first (due to Roudneff) cleverly chooses a function whose minimum is taken on a Tverberg partition.

\smallskip
{\bf Proof} by Roudneff. We assume that the points of $X$ are in general position (the coordinates are algebraically independent, say). Assume $\p=\{X_1,\ldots,X_r\}$ is an $r$-partition of $X$ with $1\le |X_j|\le d+1$ and define the function
\[
f(x,\p)=\sum_{j=1}^r \dist^2(x,\conv X_j).
\]
Here $\dist$ is the distance given by the Euclidean norm, which is denoted by $\|\cdot\|$. 
For fixed $\p$ the function $f$ is convex on $\rr^d$.  It tends to infinity as $\|x\| \to \infty$ so it attains its minimum. Then there is a partition, say $\p$, where the minimum of the function $f(x,\p)$ is the smallest; let it be $\mu$. We are going to show that $\mu=0$, which clearly suffices. Assume on the contrary, that $\mu >0$ and is reached at $z \in \rr^d$. Denote by $y_j$ the (unique) point in $\conv X_j$ with $\dist(z,\conv X_j)=\|z-y_j\|$. The function $x \mapsto  \sum_1^r \|x-y_j\|^2$ takes its minimum also at $x=z$ so its gradient at $x=z$ is zero: $\sum_1^r (z-y_j)=0$. Note that $z=y_j$ is possible but cannot hold for all $j$ since $\mu>0$.

Define $Y_j\subset X_j$ for $j=1,\ldots,r$ via $y_j \in \relint \conv Y_j$. We {\bf claim} that $\bigcap _1^r \aff Y_j=\emptyset$. Otherwise there is a point $v \in \bigcap _1^r \aff Y_j$.  Let $\langle \cdot, \cdot \rangle$ denote the standard scalar product, so $\langle x,x \rangle=\|x\|^2$, for instance. Then $\langle z-v,z-y_j\rangle>0$ if $y_i \ne z$ (because $y_j$ is the closest point to $z$ in $\conv Y_j$) and  $\langle z-v, z-y_j\rangle =0$ if $y_j = z$. Summing these inequalities and equalities gives $\langle z-v,\sum_1^r (z-y_j)\rangle >0$,  contradicting $\sum_1^r (z-y_j)=0$.

The dimension counting argument mentioned in the introduction shows now that $\sum_1^r |Y_j|\le (r-1)(d+1)$ so one point of $X$, say $x$, is not used in any $Y_j$. This is the point where the general position of $X$ is used. We can decrease the value $\mu$ if $\langle x-y_j,z-y_j \rangle>0$ for some $j$ with $y_j\ne z$ because by adding $x$ to $Y_j$ there appears a point on the segment $[x,y_j]\subset \conv (Y_j\cup \{x\})$ that is closer to $z$ than $y_j$. Thus $\langle x-y_j, z-y_j\rangle \le 0$ must hold for every $j$. Summing these inequalities gives
\begin{eqnarray*}
0 &\ge& \sum_1^r \langle x-y_j, z-y_j\rangle = \sum_1^r \Big\langle(x-z)+(z-y_j),z-y_j\Big\rangle \\
   &=& \Big\langle x-z, \sum_1^r (z-y_j)\Big\rangle+ \sum_1^r \langle z-y_j, z-y_j \rangle =0+\mu >0,
\end{eqnarray*}
a contradiction.
\qed

\smallskip
{\bf Proof} by Sarkaria. This proof has two ingredients. One is the so-called Colorful Carath\'eodory theorem of the first author \cite{Barany:1982va}. Carath\'eodory's classical theorem~\cite{Carath1907} says in essence that being in the convex hull has a very finite reason. Precisely, if $A \subset \rr^d$ and $a \in \conv A$, then $a \in \conv B$ for some $B \subset A$ with $|B| \le d+1$.  In the colorful version there are $d+1$ sets or ``colors'' $A_1,\dots,A_{d+1} \subset \rr^d$ and $a \in \bigcap_{i=1}^{d+1} \conv A_i$. A {\sl transversal} of the sets $A_1,\dots,A_{d+1}$ is simply a set with a point $a_i \in A_i$ for every $i$.

\begin{theorem}\label{th:Car} Assume $A_1,\dots,A_{d+1} \subset \rr^d$ and $a \in \bigcap_{i=1}^{d+1} \conv A_i$. Then there is a transversal $\{a_i \in A_i: i \in [d+1]\}$, such that $a \in \conv\{a_1,\dots,a_{d+1}\}$.
\end{theorem}

The colorful version contains the original one: simply take $A_i=A$ for every $i$.

\begin{figure}
\centerline{\includegraphics[scale=1]{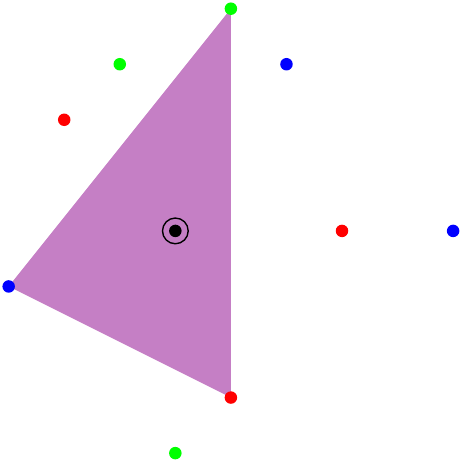}}
\caption{The colorful Carath\'edory theorem in dimension two.  Every color class contains the origin in its convex hull.  The figure shows a colorful transversal that preserves this property.}
\end{figure}

The second ingredient is Sarkaria's tensor trick~\cite{Sarkaria:1992vt}. We explain it in the form given in  \cite{BaranyOnn}. It begins with an artificial tool: choose vectors $v_1,\dots,v_r \in \rr^{r-1}$ so that their unique (up to a multiplier) linear dependence is $v_1+\dots+v_r=0$. Now let $X=\{x_0,x_1,\ldots,x_n\}$ be the set of $(r-1)(d+1)+1$ points given in Tverberg's theorem, so $n=(r-1)(d+1)$. With $x_i$ and $v_j$ we associate the tensor
\[
\overline x_{i,j}= v_j \otimes (x_i,1) \in \rr^n,
\]
the tensor $\overline x_{i,j}$ can be thought of as an $(r-1)\times (d+1)$ matrix as well. Note that we moved to the $n$-dimensional space because $\overline x_{i,j}\in \rr^n$, while the original points $x_i$ are in $\rr^d$. Observe that the origin is in the convex hull of the set
\[
A_i=\{\overline x_{i,1},\overline x_{i,2},\ldots,\overline x_{i,r}\}
\]
for every $i$. The Colorful Carath\'eodory theorem applies now in $\rr^n$ and gives, for every $x_i$, a tensor $\overline x_{i,j(i)}$ with $0\in \conv \{\overline x_{0,j(0)},\overline x_{1,j(1)},\ldots,\overline x_{n,j(n)}\}$. Thus $0 \in \rr^n$ can be written as a convex combination of the tensors $\overline x_{i,j(i)}$:
\begin{eqnarray*}
0&=&\sum_{i=0}^n \al_i \overline x_{i,j(i)}= \sum_{i=0}^n \al_i v_{j(i)} \otimes (x_i,1)\\
  &=& \sum_{j=1}^r  v_j \otimes \left( \sum_{i:j=j(i)}\al_i(x_i,1)\right) = \sum_{j=1}^r  v_j \otimes \left( \sum_{x_i\in X_j}\al_i(x_i,1)\right),
\end{eqnarray*}
where $X_j:=\{x_i\in X: j(i)=j\}$. These sets form a partition of $X$ into $r$ parts. There is a vector $u \in \rr^{r-1}$ orthogonal to $v_3,\ldots,v_r$ such that $\langle u, v_1\rangle =1$. Then $\langle u,v_2\rangle =-1$ because of the condition  $v_1+\dots+v_r=0$. Multiplying the last equation by $u$ from the left gives $\sum_{x_i\in X_1}\al_i(x_i,1)=\sum_{x_i\in X_2}\al_i(x_i,1)$. It follows then that
\[
\sum_{x_i\in X_1}\al_i(x_i,1)=\sum_{x_i\in X_2}\al_i(x_i,1)=\ldots=\sum_{x_i\in X_r}\al_i(x_i,1).
\]
Reading the last coordinate here shows that $\al:=\sum_{x_i\in X_1}\al_i=\sum_{x_i\in X_2}\al_i=\ldots=\sum_{x_i\in X_r}\al_i >0$. (Actually $\al=1/r$.) Then
\[
p:=\frac 1{\al}\sum_{x_i\in X_1}\al_ix_i=\frac 1{\al}\sum_{x_i\in X_2}\al_ix_i=\ldots=\frac 1{\al}\sum_{x_i\in X_r}\al_ix_i
\]
is a point in the convex hull of each $X_j$: $X_1,\ldots,X_r$ is the required partition.
\qed

There is more to Sarkaria's method than just this proof. To see this let $X_1,\ldots,X_r$ be finite (or compact) sets in $\rr^d$. What condition guarantees that $\bigcap_1^r \conv X_j = \emptyset$? There is a classical necessary and sufficient condition:

\begin{theorem}\label{th:halfspace} Under the above conditions,  $\bigcap_1^r \conv X_j=\emptyset$ if and only if there are closed halfspaces $D_1,\dots,D_r$ with $X_j \subset D_j$ for every $j\in [r]$ such that  $\bigcap_1^r D_j=\emptyset$.
\end{theorem}

The {\bf proof} is easy. One direction is trivial. In the other direction the case $r=2$ is just the separation theorem for convex sets, and induction on $r$ works for $r>2$.

Here comes another necessary and sufficient condition from Arocha et al~\cite{Arocha:2009ft}. First define $X=\bigcup_1^r X_j$, here either $X$ is a multiset or we assume that the sets $X_j$ are disjoint.  For $x \in X$ denote, as before,
\[
\overline{x} = v_j \otimes (x,1) \mbox{ if } x\in X_j \mbox{ and set } \overline{X} = \{\overline{x} : x \in X\}.
\]
Here the vectors $v_j \in \rr^{r-1}$ are the same as before.

\begin{theorem}\label{th:sark} Under the above conditions, $\bigcap_1^r \conv X_j \ne \emptyset$ if and only if $0 \in \conv \overline X$.
\end{theorem}

The {\bf proof} is essentially the same as above, starting with the convex combination of the vectors in $\overline{X}$ representing the origin:
\begin{eqnarray*}
0&=&\sum_{x \in X} \al(x)\overline{x}=\sum_{j=1}^r\sum_{x \in X_j}\al(x)v_j\otimes (x,1)\\
  &=&\sum_{j=1}^rv_j\otimes \sum_{x \in X_j}\al(x)(x,1).
\end{eqnarray*}
After this factorization the arguments are analogous to the previous proof.
\section{Topological versions}\label{section-topological}

\begin{figure}[h!]
\centerline{\includegraphics[scale=0.7]{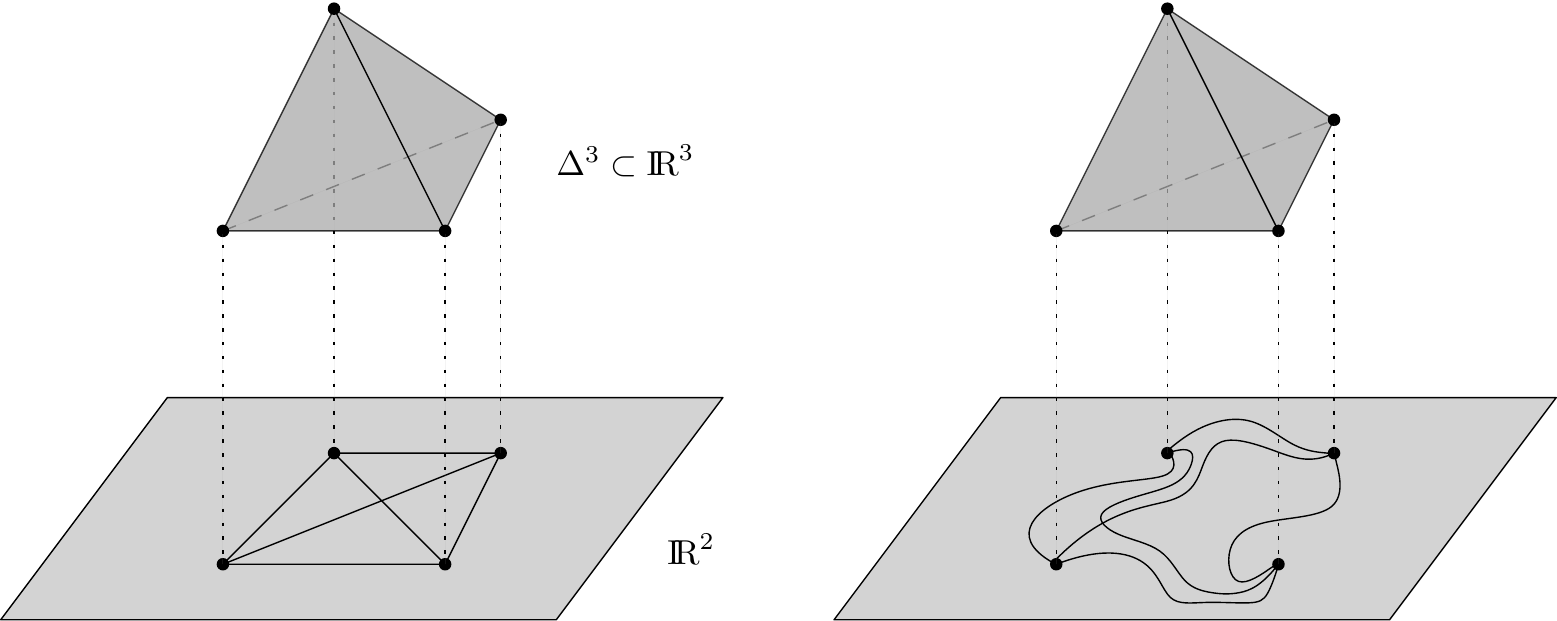}}
\caption{Radon's theorem, affine and topological versions in $\rr^2$}
\label{figtop}
\end{figure}

We start with a different formulation of Radon's theorem. Given a set $X$ of $d+2$ points in $\rr^d$, there is a $(d+1)$-dimensional simplex $\tri^{d+1}$ with vertex set $V$ and an affine map $f:\rr^{d+1} \to \rr^d$ such that $f(V)=X$. Proper faces of the simplex are mapped to the convex hull of the corresponding points of $X$. So Radon's theorem says, in this setting, that {\sl there are disjoint (proper) faces of $\tri^{d+1}$ whose $f$-images intersect}, see Figure \ref{figtop}. What happens if $f: \tri^{d+1} \to \rr^d$ is not affine, but only continuous? The answer is the following theorem of Bajm\'oczy and B\'ar\'any from 1979~\cite{Bajmoczy:1979bj}, where $\skel_k \tri^{d+1}$ denotes the $k$-dimensional skeleton of $\tri^{d+1}$.

\begin{theorem}[Topological Radon]\label{th:tRadon} If $f:\skel_d \tri^{d+1} \to \rr^d$ is continuous, then the simplex has two disjoint faces $\si_1,\si_2$ with $f(\si_1)\cap f(\si_2)\ne \emptyset$.
\end{theorem}

In other words, the $d$-skeleton of the $(d+1)$-simplex cannot be embedded in $d$-space without mapping two points from disjoint faces to the same point in $\rr^d$. Actually, this holds for any $(d+1)$-dimensional polytope, not only for the simplex. In this form the result is used (and proved in a slightly more general form) by Lov\'asz and Schrijver~\cite{Lovasz:1989} in connection with the Colin de Verdi\'ere number of graphs. The proof of Theorem~\ref{th:tRadon} uses the Borsuk-Ulam theorem.

The famous non-embeddability theorem of Van Kampen and Flores~\cite{van1933komplexe, Flores1932cx} says that {\sl the $d$-skeleton of the $(2d+2)$-dimensional simplex cannot be embedded in $\rr^{2d}$}. The particular case of $d=1$ is half of Kuratowski's theorem on planar graph: the complete graph $K_5$ on five vertices is not planar. Sarkaria~\cite{Sarkaria:1991ug} realized in 1991 that there is some connection between the topological Radon theorem and the Van Kampen and Flores theorem. Recently it has been shown by Blagojevi\'c, Frick and Ziegler~\cite{Blagojevic:2014js} that the topological Radon theorem implies Van Kampen-Flores.  The same implication is mentioned (somewhat implicitly) in Gromov~\cite{Gromov:2010eb} as well. The proof is by the {\sl constraint method}, a powerful new technique that has several further implications. Here is how it goes in the given case.

{\bf Proof.} Assume that there is a map $f: \skel_d \tri^{2d+2} \to \rr^{2d}$ that sends any two points from disjoint faces of  $\skel_d \tri^{2d+2}$ to distinct points in $\rr^d$. Extend this map to the $2d+1$ skeleton of $\tri^{2d+2}$ continuously (but otherwise arbitrarily) and define a new map
\[
g: \skel_{2d+1} \tri^{2d+2} \to \rr^{2d+1}
\]
where the first $2d$ coordinates of $g(x)$ coincide with those of $f(x)$ and the last coordinate of $g(x)$ is simply the distance of $x$ from the $\skel_d \tri^{2d+2}$. Since $g$ is continuous, the topological Radon shows now that for some two points, say $x_1$ and $x_2$ from pairwise disjoint faces of $\tri^{2d+2}$, $g(x_1)=g(x_2)$.
So $f(x_1)=f(x_2)$ and $\dist(x_1,\skel_d \tri^{2d+2})=\dist(x_2,\skel_d \tri^{2d+2})$. But as $x_1$ and $x_2$ belong to disjoint faces, one of these faces is of dimension at most $d$, so the last components of both $g(x_1)$ and $g(x_2)$ are equal to zero, that is, both $x_1,x_2 \in \skel_d \tri^{2d+2}$.
\qed

Of course Tverberg's theorem can be reformulated the same way: if $f: \skel_d \tri^n \to \rr^d$ is an affine map and $n=(r-1)(d+1)$, then there are disjoint faces $F_1,\ldots,F_r$ of $\tri^n$ such that $\bigcap _1^r f(F_j)\ne \emptyset$. This statement is equivalent to Tverberg's theorem. The continuous version, a question of the first author from 1976, had been a conjecture for almost 40 years. On the positive side, the following is known.

\begin{theorem}[Topological Tverberg]\label{th:tTverb} If $f: \skel_d \tri^n \to \rr^d$ is a continuous map,  $n=(r-1)(d+1)$, and $r$ is a prime power, then there are disjoint faces $F_1,\ldots,F_r$ of $\tri^n$ such that $\bigcap _1^r f(F_j)\ne \emptyset$.
\end{theorem}

The case when $r$ is prime was proved by B\'ar\'any, Shlosman, Sz\H ucs~\cite{Barany:1981vh}  in 1981, and the prime power case by \"Ozaydin~\cite{Oza87} in 1987 in an unpublished yet influential paper, see also \cite{Volovikov:1996up}. We now give the sketch of the proof for the case when $r$ is prime. 

{\bf Proof.} Assume  $f: \skel_d \tri^n \to \rr^d$ is a counterexample. Consider the $r$-fold deleted product $D(n,r)$ of $\tri^n$, that is, the set of $r$-tuples $(x_1,\ldots,x_r)$ where the points $x_j$ come from disjoint faces of $\tri^n$. Then the map $F : D(n,r) \to \rr^{dr}$ defined by
\[
F(x_1,\ldots,x_r)=(f(x_1),\ldots,f(x_r))
\]
avoids the diagonal $\{(x,\ldots,x)\in \rr^{dr}: x\in \rr^d\}$. Note that the cyclic group  $\Z_r$ acts on the spaces $D(n,r)$ and $\rr^{dr}$: its generator $\omega$ maps $(x_1,\ldots,x_r)\in D(n,r)$ to $(x_2,\ldots,x_r,x_1)$ and $(z_1,\ldots,z_r)\in \rr^{dr}$ to $(z_2,\ldots,z_r,z_1)$. Moreover, $F$ is $\Z_r$-equivariant:
\[
F(\omega(x_1,\ldots,x_r))=\omega (F(x_1,\ldots,x_r)).
\]

Actually, the symmetric group $S_r$ on $r$ elements acts on $D(n,r)$ and $\rr^{dr}$ equivariantly as well but, for this proof, the action of its subgroup $\Z_r$ suffices. The orthogonal complement of this diagonal is $W(n,r)=\{(x_1,\ldots,x_r)\in \rr^{dr}: x_1+ \ldots +x_r=0\}$ which is in fact isomorphic to $\rr^{d(r-1)}$; its unit sphere is $S(W(n,r))$. Consider the chain of maps
\[
D(n,r) \to \rr^{dr}\setminus \operatorname{diagonal} \to W(n,r)\setminus \{0\} \to S(W(n,r))
\]
where $\rr^{dr} \to W(n,r)$ is the orthogonal projection onto the subspace $W(n,r)$ and the map $W(n,r) \to S(W(n,r))$ sends $x \in W(n,r)$, $x\ne 0$ to $x/\|x\| \in S(W(n,r))$. The composition is a map $G: D(n,r) \to S(W(n,r))$ which is again $\Z_r$-equivariant. In addition, the action of $\Z_r$ is {\sl free} on both $D(n,r)$ and $S(W(n,r))$, meaning that the orbit of any point in $D(n,r)$ and $S(W(n,r))$ consists of $r$ distinct points. This is because $r$ is prime. In this case Dold's theorem, an extension of the Borsuk-Ulam theorem, applies: there is no $\Z_r$-equivariant map from an $(n-r)$-connected space to an $(n-r)$-dimensional space provided the action is free on both spaces \cite{Dold:1983wr}. Here $S(W(n,r))$ is $(n-r)$-dimensional trivially, and, as shown in \cite{Barany:1981vh}, $D(n,r)$ is $(n-r)$-connected.\qed

\medskip
This proof is a typical example of the {\sl configuration space - test map scheme} (consult \cite{Matousek:2002td} and the references therein for more on this method). When this is applied for the prime power $r=p^k$ case of the topological Tverberg theorem, the map  $G: D(n,r) \to S(W(n,r))$ is equivariant with respect to the abelian group $\left(\Z_{p}\right)^k$ but the action is not free. What \"Ozaydin~\cite{Oza87} observes is that it is fixed point free (i.e., no point is fixed by all group elements) and so some algebraic topology machinery still works and excludes the existence of such a map. \"Ozaydin goes one step further and shows that, if $r$ is not a prime power, then there is a map  $D(n,r) \to S(W(n,r))$ which is equivariant under the symmetric group $S_r$. Consequently, the configuration space - test map scheme fails badly here. So what comes next? Is there a topological Tverberg Theorem for non-prime $r$? This had been ``one of the most
challenging problems in this field'' according to Matou\v{s}ek, and ``a holy grail of topological combinatorics'' according to Kalai. This question had remained open for almost forty years.

In 2010, in a groundbreaking paper, Gromov~\cite{Gromov:2010eb} states that ``The topological Tverberg theorem, whenever available, implies the (generalized) Van Kampen-Flores theorem''. This implication holds for any $r$, prime or not. Gromov also gives a proof (or rather a sketch of proof) in three lines. A detailed proof can be found in ~\cite{Blagojevic:2014js}. Surprisingly, this remark of Gromov went completely unnoticed.

The generalized Van Kampen-Flores theorem is due to Sarkaria~\cite{Sarkaria:1991ug} when $r$ is prime and to Volovikov~\cite{Volovikov:1996up} when $r$ is prime power. It says the following.

\begin{theorem}[Generalized Van Kampen-Flores]\label{genVKamFl} Let $d \ge 1$ be an integer, let $r$ be a prime power, let $k\ge (r-1)d/r$ be an integer, $N = (d + 2)(r-1)$, and let $f: \tri^N \to \rr^d$ be a continuous map. Then there exist $r$ pairwise disjoint faces $\si_1,\ldots,\si_r$ in the $k$-skeleton of the simplex $\tri^N$ whose $f$-images overlap: $f(\si_1) \cap \ldots \cap f(\si_r) \ne \emptyset$.
\end{theorem}

The proof, rediscovered by Blagojevi\'c, Frick, and Ziegler~\cite{Blagojevic:2014js} is almost identical to the previous proof for the case $r=2$. It has two ingredients: one is the topological Tverberg theorem, the other is the constraint method (or the pigeonhole principle). The proof also works when $r$ is not a prime power and shows that if the generalized Van Kampen-Flores theorem fails, then so does the topological Tverberg.

Unaware of Gromov's remark connecting the topological Tverberg and  the generalized Van Kampen-Flores theorems, Mabillard and Wagner started working on extending the Whitney trick~\cite{whitney1944} to an $r$-fold Whitney trick. Their hope was that the method, when combined with \"Ozaydin's example, would give a counterexample to the topological Tverberg conjecture in the non-prime power case. What they proved is the following remarkable result~\cite{MabillardWagner2015}.

\begin{theorem} Let $K$ be an $(r-1)\ell$-dimensional simplicial complex where $r\ge2,\;\ell\ge 3$ are integers, and let $D(K,r)$ denote the $r$-fold deleted product of $K$. Then the following two statements
are equivalent:
\begin{itemize}
\item there exists an $S_r$-equivariant map $D(K,r) \to S(W(r\ell,r))$,
\item there exists a continuous map   $f:K \to \rr^{r\ell}$ such that
the $f$-images of any $r$ disjoint faces of $K$ have no point in common.
\end{itemize}
\end{theorem}

Mabillard and Wagner almost succeeded in finding a counterexample to the topological Tverberg conjecture: what was missing was an example where the generalized Van Kampen-Flores theorem fails. It was Florian Frick~\cite{Frick2015} who realized that the above theorem and \"Ozaydin's example combined with the constraint method (or Gromov's remark), gives a counterexample for every non-prime power $r$. A more detailed description, with further applications of the constraint method is presented in \cite{Blagojevic:2017js}.

The specific example in \cite{Frick2015} is with $r=6$ and $d=19$, so there is a continuous map from the 19-skeleton of $\tri^{100}$ to $\rr^{19}$ such that the images of any 6 disjoint faces have no point in common. Subsequently this was further improved, by Avvakumov, Mabillard, Skopenkov, and Wagner~\cite{AMSW15eliminating} to a map $\tri^{65} \to \rr^{12}$ with the same property. Moreover, Mabillard and Wagner~\cite{MabillardWagner2015} came up with another counterexample without using Gromov's or Blagojevi\'c et al. reduction. 

There is hope for positive topological results related to Tverberg's theorem if $r$ is not a prime power.  Even though Tverberg partitions may not exist, strong intersection properties of the images of disjoint faces of $\tri^n$ under a map $f:\tri^n \to \rr^d$ can be obtained \cite{Simon16}.  If we are allowed to use more points, a topological version of Birch's theorem is still open.

\begin{problem}
Decide if the following statement is true.  If $f: \skel_d \tri^n \to \rr^d$ is a continuous map,  $n=r(d+1)-1$, then there are disjoint faces $\si_1,\ldots,\si_r$ of $\tri^n$ such that $\bigcap_{j=1}^r f(\si_j)\ne \emptyset$.
\end{problem}

This was first presented as Conjecture 5.5 in \cite{Blagojevic:2014js}, where it is also explained that $n=(r-1)d-1$ is the smallest value where the conjecture could conceivably be true.

\section{Colorful versions}\label{section-colored}

One intriguing family of variations of Tverberg-type theorems is the colorful versions of Tverberg's theorem.  The goal is to restrict to partitions of a set $X$ of points where some pairs of points are required to be in different parts.  This is usually achieved by coloring the points with few colors, and asking that no part in the partition has more than one point of any color. Motivation came from the halving plane problem as  explained in Section \ref{section-applications}. The main open problem of this kind is the following variant of a conjecture by B\'ar\'any and Larman \cite{Barany:1992tx}.

\begin{conjecture}[Colored Tverberg theorem]\label{conjecture-colored92}
Let $r,d$ be positive integers.  Let $t=t(d,r)$ be the smallest positive integer, if exists, such that for any $d+1$ sets $F_1, F_2, \ldots, F_{d+1}$ of $t$ points each in $\rr^d$, considered as color classes, there are $r$ disjoint sets $X_1, \ldots, X_r$ of $X = \cup_{i=1}^{d+1} F_i$ such that
\begin{itemize}
	\item each $X_j$ has exactly one point of each $F_i$ and
	\item the convex hulls of the sets $X_j$ intersect.
\end{itemize}
Then, for any $r,d$, the number $t(d,r)$ exists and is equal to $r$.
\end{conjecture}

\begin{figure}
\centerline{\includegraphics[scale=1]{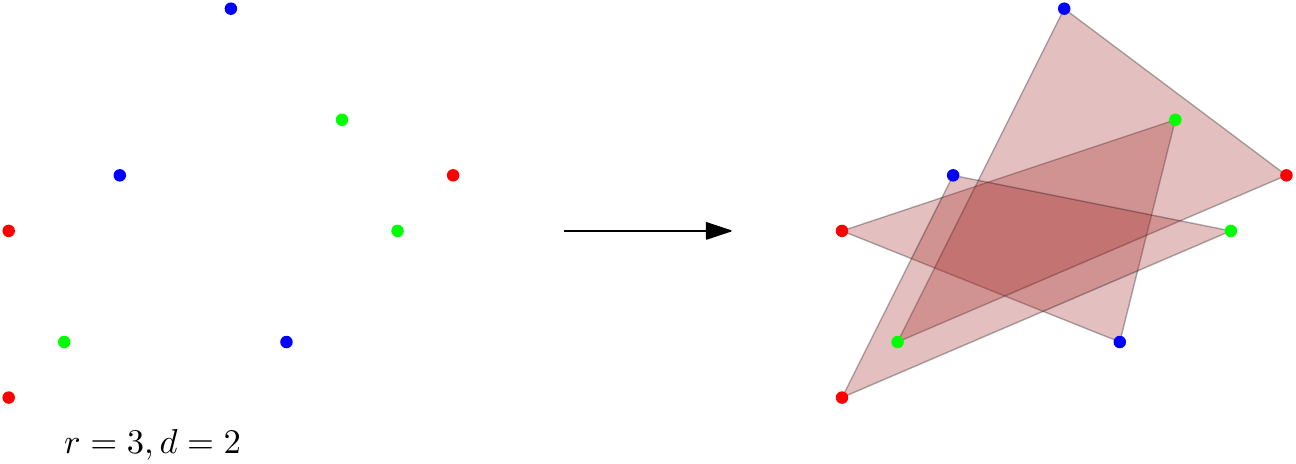}}
\caption{A colorful Tverberg partition.}
\end{figure}

In the original conjecture in \cite{Barany:1992tx} each color class is of size at least $r$, and the question is whether there is an integer $n(d,r)$ with the following property. If the union of the color classes is of size $n(d,r)$, then there are disjoint sets $X_1, \ldots, X_r \subset\cup_{i=1}^{d+1} F_i$
satisfying the two conditions above. This conjecture is still open in general, but clearly $n(d,r)=r(d+1)$ if $t(d,r)=r$.

A partition as above is called a {\it colorful Tverberg partition}. The first result of this kind was obtained by B\'ar\'any, F\"uredi, and Lov\'asz, showing that $t(2,3) \le 7$ \cite{Barany:1990wa}.  In the paper containing Conjecture \ref{conjecture-colored92}, B\'ar\'any and Larman showed that it is true for $d=2$ and any $r$.  Lov\'asz proved the case $r=2$ and any $d$, also known as the colorful Radon theorem, using the Borsuk-Ulam theorem.  His proof appears in \cite{Barany:1992tx}.  Here we include a linear-algebraic proof of the colorful Radon theorem, from \cite{soberon2015equal}.

\smallskip
{\bf Proof.} If we are given $F_1, \ldots, F_{d+1}$ pairs of points in $\rr^d$, we can name their elements arbitrarily $F_i = \{x_i, y_i\}$ for each $i$.  Consider the $d+1$ vectors of the form $x_i-y_i$.  Since we have more than the dimension, they must have a non-trivial linear dependence,
\[
\sum_{i=1}^{d+1}\alpha_i (x_i - y_i) = 0.
\]
If there is any $\alpha_i < 0$, we can swap the names of $x_i$ and $y_i$ and the sign of $\alpha_i$ without breaking the linear dependence.  Once all signs are non-negative, we can assume by scaling that their sum is $1$, as they were not all zero.  A simple manipulation of the linear dependence gives
\[
\sum_{i=1}^{d+1}\alpha_i x_i = \sum_{i=1}^{d+1}\alpha_i y_i.
\]
Thus, the partition $A=\{x_1, \ldots, x_{d+1}\}$, $B= \{y_1, \ldots, y_{d+1}\}$ satisfies the requirements.
\qed

The existence of $t(d,r)$ was first settled by \v{Zivaljevi\'c} and Vre\'cica \cite{Zivaljevic:1992vo}, showing that $t(d,r) \le 2r-1$ if $r$ is a prime number, which implies $t(d,r) \le 4r-3$ for all $r$.  The proof is topological, and extends to the topological version of the colorful Tverberg theorem.

One thing to be noted about Conjecture \ref{conjecture-colored92} is that it does not imply Tverberg's theorem directly.  For colorful versions of other classic results, such as Carath\'eodory's theorem or Helly's theorem, when the color classes are equal we recover the original result \cite{Barany:1982va}.  However, there is a theorem by Blagojevi\'{c}, Matschke and Ziegler, also called the \textit{optimal colorful Tverberg}, which generalizes both Tverberg's theorem and the B\'ar\'any-Larman conjecture.  They proved the following theorem and its topological analogue.  

\begin{theorem}[Optimal colorful Tverberg \cite{blago15, blago-optimal-2011}]\label{theorem-optimalcolored}
Let $r$ be a prime number, $n=(r-1)(d+1)+1$, and $X$ be a set of $n$ points in $\rr^d$.  Suppose that they are colored with $c$ colors so that each color appears at most $r-1$ times.  Then, there is a partition of $X$ into $r$ parts $X_1, \ldots, X_r$ so that each $X_j$ has at most one point of each color and their convex hulls intersect.
\end{theorem}

This implies Conjecture \ref{conjecture-colored92} when $r+1$ is a prime number in the following way.  Given $F_1, \ldots, F_{d+1}$ sets of $r$ points each in $\rr^d$, add an extra point $p_0$ of a new color.  Then, we can apply Theorem \ref{theorem-optimalcolored} to the total set of $r(d+1)+1$ points.  This gives us a partition into $r+1$ parts, and we can simply drop the part containing $p_0$ and redistribute its points to have a partition as in Conjecture \ref{conjecture-colored92}.

\begin{problem}
Does Theorem \ref{theorem-optimalcolored} hold for all $r>1$?
\end{problem}

\begin{problem}
Is there a non-topological proof of the colorful Tverberg theorem for $r \ge 3$?
\end{problem}

The constraint method by Blagojevi\'{c}, Ziegler, and Frick \cite{Blagojevic:2014js}, mentioned in Section \ref{section-topological} also gives colorful Tverberg's results.  We showcase here how it implies the bound $t(d,r) \le 2r-1$ when $r$ is a prime power.  The reader may notice that both the statement and the proof carry on through the topological setting.

\begin{theorem}
	Let $r$ be a prime power and $n = (2r-2)(d+1)+1$.  Suppose that the vertices of $\Delta^{n-1}$ are colored with $d+1$ colors, each of which appears at most $2r-1$ times.  Then, for any continuous function $f: \Delta^{n-1} \to \rr^d$, there are points $x_1, \ldots, x_r \in \Delta^{n-1}$ in pairwise vertex-disjoint faces such that each $x_i$ is contained in a face that has at most one point of each color and $f(x_1) = \ldots = f(x_{r})$.
\end{theorem}

\smallskip
{\bf Proof.} For each color $i$, let $M_i$ be the simplicial complex of faces of $\Delta^{n-1}$ with at most one vertex of color $i$.  We can then define
\begin{align*}
f_i : & \Delta^{n-1} \to \rr \\
f_i(x) & = \operatorname{dist}(x, M_i).
\end{align*}
We can use these functions to extend $f: \Delta^{n-1} \to \rr^{d}$ to a new function $\tilde{f} = (f,f_1, \ldots, f_{d+1}): \Delta^{n-1} \to \rr^{2d+1}$.  Notice that $n = (2r-2)(d+1) + 1 = (r-1)(2d+2)+1$. Then, we can apply the colorful Tverberg theorem and find $x_1, \ldots, x_{r} \in \Delta^{n-1}$ points contained in vertex-disjoint faces such that $f(x_1) = \ldots = f(x_r)$ and $f_i (x_1) = \ldots = f_i(x_r)$ for each $i$.  However, since there are at most $2r-1$ vertices of color $i$, one of the $x_j$ must be in $M_i$.  This implies that all $x_j$ are in $M_i$, as desired.
\qed

If we have fewer than $d+1$ color classes in the colorful Tverberg theorem, we cannot guarantee the existence of a colorful Tverberg partition into $r$ parts for sufficiently large $r$. This follows simply because colorful simplices have positive co-dimension.  However, when the co-dimension is not a problem, a similar proof to the one above yields the following result and its natural topological version \cite{VZ93-new, Blagojevic:2014js}.

\begin{theorem}
	Let $r,d,c$ be positive integers such that $d > r (d+1-c)$ and $r$ is a prime power.  For any $c$ sets $F_1, \ldots, F_c$ of $2r-1$ points each in $\rr^d$, considered as color classes, we can find $r$ colorful sets $X_1, \ldots, X_r$ whose convex hulls intersect.	
	\end{theorem}

On the other hand, there are some benefits of increasing the number of color classes in Conjecture \ref{conjecture-colored92}.  In the proof we presented for colorful Radon, a careful reader may notice that we did not only find two colorful intersecting simplices, but they used the same coefficients for the convex combination that witnesses the intersection.  If we seek this for colorful partitions with $r>2$, then $(r-1)d+1$ color classes are sufficient and necessary \cite{soberon2015equal}.  The topological version of this statement also holds when $r$ is a prime power\cite{Blagojevic:2014js}.

A different way to generalize Tverberg' theorem stems from the following result by A. P\'or \cite{attila97thesis}.

\begin{theorem}\label{Kirchberger}
Given $r$ sets $X_1, \ldots, X_r \subset \rr^d$, we have that $\bigcap_{j=1}^r \conv X_j \neq \emptyset$ if and only if for every set $A \subset \bigcup_{1}^r X_j$ of at most $(r-1)(d+1)+1$ points, we have $\bigcap_{j=1}^r \conv (A \cap X_j) \neq \emptyset$.
\end{theorem}

The case $r=2$ is Kirchberger's theorem \cite{Kirch03}.  A colorful version of this result was proven in \cite{Arocha:2009ft}, which generalizes both P\'or's result and Tverberg's theorem.

\begin{theorem}
Given positive integers $r,d$, let $n=(r-1)(d+1)+1$.  We are given $n$ sets $G_1, G_2, \ldots, G_n$ which are colorful, using $r$ colors.  Let $X_1, \ldots, X_r$ be the colors classes.  Then, if every transversal $Y = \{y_1, \ldots, y_n\}$, where $y_i \in G_i$ for all $i$, satisfies that
\[
\bigcap_{j=1}^r \conv(Y \cap X_j) = \emptyset,
\]
there must be a set $G_i$ such that
\[
\bigcap_{j=1}^r \conv(G_i \cap X_j) = \emptyset.
\]
\end{theorem}

The result above implies Tverberg's theorem if each $G_i$ consists of $r$ copies of a point $a_i$, all of different colors.  It implies P\'or's result if for all $i,j$ we have $G_i \cap X_j = X_j$.

Recently, it has been observed that some classic colorful theorems in combinatorial geometry can be generalized by using matroids instead of color classes.  Examples are Kalai and Meshulam's generalization of colorful Helly \cite{Kalai:2005tb} or Holmsen's generalization of colorful Carath\'eodory \cite{holmsen2016intersection}. Such a version exists for Tverberg's theorem, as was proven by B\'ar\'any, Kalai, and Meshulam \cite{matroid-tverberg16}.

\begin{theorem}\label{theorem-tverbergmatroid}
	Let $d$ be a positive integer, $M$ be a matroid of rank $d+1$ and $b(M)$ the maximal number of pairwise disjoint bases of $M$.  Then, for any continuous map from the matroidal complex of $M$ to $\rr^d$ there are $\left\lceil \frac{\sqrt{b(M)}}{4}\right\rceil$ disjoint independent sets whose images under $f$ intersect.
\end{theorem}

Given a set $S$ whose elements are colored with $d+1$ colors, we can define the matroid $M$ on $S$ by saying that a subset is independent if it has at most one element of each color.  An application of Theorem \ref{theorem-tverbergmatroid} to $M$ yields results along the lines of Conjecture \ref{conjecture-colored92}.  However, being able to use any matroid gives much more flexibility.  More on this result can be seen in \cite{BHZ2017tverberg}.

Yet another way to impose conditions on Tverberg's theorem is using a graph.  We say that a graph $G$ on $N$ vertices is an $r$-Tverberg graph for $\rr^d$ if the following holds.  For any set of $N$ points in $\rr^d$ representing the vertices of $G$, there is a Tverberg partition of the points into $r$ parts so that each part is an independent subset of $G$.  Conjecture \ref{conjecture-colored92} can be rephrased as saying that the disjoint union of $d+1$ complete graphs $K_r$ is an $r$-Tverberg graph for $\rr^d$.  Sparser graphs than this one are known to be $r$-Tverberg graphs,  as the following result by Hell shows \cite{Hell:2008em}.

\begin{theorem}
	Let $r$ be a prime power.  Then, every graph for which each connected component either
	\begin{itemize}
		\item has cardinality smaller than $\frac{r+2}{2}$,
		\item is a complete bipartite graph $K_{1,l}$ for $l < r-1$,
		\item is a path (if $r>3$), or
		\item is a cycle (if $r>4$)
	\end{itemize}
	is an $r$-Tverberg graph for $\rr^d$.
\end{theorem}

\section{The structure of Tverberg partitions}\label{section-intersection}

Once the existence of Tverberg partitions has been established, the next step is to have a better understanding of their structure. One way is to relax the conditions on the partition, such as asking for the convex hulls of the parts to have a transversal low-dimensional affine subspace (as opposed to a point in common) or for the parts to have pairwise intersection.  Another is to strengthen the conclusion of the theorem, such as guaranteeing many Tverberg partitions, seeking partitions which are resistant to changes in the point set, or determining the dimension of the set of points which witness the intersection of a Tverberg partition.

\subsection{Sierksma's conjecture}

One of the most notable open problems around Tverberg's theorem is to give a lower bound for the number of Tverberg partitions we can find in any set of $(r-1)(d+1)+1$ points.  Tverberg's theorem shows that at least one partition always exists, but in general there is hope for much more.  This was formalized by Sierskma with his now famous conjecture \cite{sierksma1979convexity}.

\begin{conjecture}
	Every set of $(r-1)(d+1)+1$ points in $\rr^d$ has at least $(r-1)!^d$ different Tverberg partitions.
\end{conjecture}

It is also known as ``the Dutch cheese conjecture'' since Sierksma promised a Dutch cheese as a prize for a solution.  The number $(r-1)!^d$ cannot be improved.  A simple example is to take the vertices of a simplex, and cluster $r-1$ points near each vertex and one final vertex in the barycenter of the simplex.  The number of Tverberg partitions can be easily counted to be $(r-1)!^d$. A large (and very different) family of examples exhibiting this bound have been constructed by White ~\cite{Whi17}. Actually White answers the following question of Perles. 

Every Tverberg partition $X_1,\ldots,X_r$ (of a set $X \subset \rr^d$ with $n=(r-1)(d+1)+1$ elements) defines a partition of $[n]$ into $r$ integers $k_1,\ldots,k_r$ where  $k_i=|X_i|$ for all $i$. Of course $k_i \in  [d+1]$. Call this partition of $[n]$ the {\it signature} of this Tverberg partition of $X$. Perles asked whether, given such a partition of $[n]$, is there a set $X \in \rr^d$ of $n$ elements such that the signature of every Tverberg partition of $X$ has the given partition of $[n]$. This was answered in the affirmative by the following interesting theorem by  White.

\begin{theorem}\label{theorem-white} Assume $d\ge 1, r\ge 2$ and $n=(r-1)(d+1)+1$. Given integers $k_1,\ldots,k_r$  with $k_i \in [d+1]$ for every $i \in [r]$ and $k_1+\ldots+k_r=n$, there is a set $X \in \R^d$ such that the signature of every Tverberg partition of $X$ coincides with the multiset $\{k_1,\ldots,k_r\}$.
\end{theorem}

It is not hard to see that the number of these Tverberg partitions is  $(r-1)!^d$. There are further  families of examples achieving this lower bound in Sierksma's conjecture.  This is explained in detail in Section \ref{section-universal}.

It is also an interesting question if Sierksma's bound holds for the topological versions when $r$ is a prime power.  As for lower bounds, Vu\v{c}i\'c and \v{Z}ivaljevi\'c proved by topological methods that one can always find $\frac{1}{(r-1)!}(r/2)^{(r-1)(d+1)/2}$ Tverberg partitions if $r$ is a prime number \cite{Vucic:1993be}, which extends to the topological version of the problem.  In rough terms this is the square root of the lower bound conjectured by Sierksma.  This was extended to prime powers by Hell \cite{Hell:2007tp}.  The only non-trivial case of the conjecture which has been verified is $d=2, r=3$ by Hell \cite{Hell:2008em}.

If $r$ is not a prime power, then we must restrict ourselves to the affine version of the problem.  In this case, the best bound is that $(r-d)!$ Tverberg partitions exist \cite{Hell08Birch}.  Bounds for the number of partitions in the colorful case (namely, for instances of Theorem \ref{theorem-optimalcolored}) and for general Birch partitions (when we use $r(d+1)$ points instead of $(r-1)(d+1)+1$) can be found in \cite{Hell14colored, Hell08Birch}.

\subsection{The Tverberg-Vre\'cica conjecture}

As mentioned in Section \ref{section-introduction}, one motivation for Tverberg's theorem is the centerpoint theorem, where one finds a point that is ``very deep'' in some finite set $X \subset \rr^d$.  One of the ways to generalize this theorem is to get a version that works simultaneously for many point sets.  This is shown in the following result, proven independently in \cite{Dolnikov:1992ut, Zivaljevic:tw}.

\begin{theorem}\label{theorem-dolnikovtransversal}
Let $0 \le k \le d-1$ be integers.  Given $k+1$ finite sets $X^1, \ldots, X^{k+1}$ of points in $\rr^d$, there is a $k$-dimensional affine subspace $L$ such that any closed half-space containing $L$ also has at least
\[
\frac{|X^m|}{d-k+1}
\]
points of $X^m$, for every $m \in [k+1]$.
\end{theorem}

Note that the case $k=0$ is the centerpoint theorem. The other end, the case $k=d-1$ is the discrete version of the classic ham-sandwich theorem. The latter is a consequence of the Borsuk-Ulam theorem and says the following. Given $d$ nice probability measures $\mu_1,\ldots,\mu_d$ in $R^d$, there is a hyperplane that splits the space into two halfspaces $H^+$ and $H^-$ so that $\mu_i(H^+)=\mu(H^-)=1/2$ for every $i \in [d]$. (A measure $\mu$ on $\rr^d$ is {\sl nice} if $\mu(h)=0$ for every hyperplane $h$.) 

Just like Tverberg's theorem is a discrete version of the centerpoint theorem, one may wonder if there is a discrete analogue of the theorem above.  This was conjectured by Tverberg and Vre\'cica \cite{TV93}.

\begin{conjecture}
Let $0 \le k \le d$ be integers.  Suppose that we are given integers $r_1, r_2, \ldots, r_{k+1}$ and sets $X^1, \ldots, X^{k+1}$ of points of $\rr^d$.  If for each $m \in [k+1]$ we have $|X^m| = (r_m - 1)(d-k+1)+1$, then we can partition each $X^m$ into $r_m$ parts $X^m_1, \ldots X^m_{r_m}$ in such a way that there is a $k$-dimensional affine subspace that intersects $\conv(X^m_j)$ for all $m \in [k+1]$ and $j \in [r_m]$.
\end{conjecture}

The case $k=0$ is Tverberg's theorem, and $k=d$ follows from taking $L=\rr^d$.  Tverberg and Vre\'cica proved the case $k=d-1$ in \cite{TV93} and a slightly weaker form of the case $k=d-2$.  The conjecture has also been verified by Karasev when all $r_m$ are powers of the same prime $p$ and $p(d-k)$ is even \cite{Kar07}, which extends two prior results \cite{Ziv99, Vre03}.  These generalizations work in the topological version of the conjecture.

There is a colorful version of the Tverberg-Vre\'{c}ica conjecture.  To see this, we consider the elements of each $X^m$ to be colorful.  We ask for the partition of each $X^m$ to satisfy that no two points of the same color are in the same part.  In the special case when $r_1 = \ldots = r_m=p$ for some prime number $p$, either $p(d-k)$ is even or $k=0$, and each $X^m$ is a colorful set such that no color class has more than $p-1$ points, Blagojevi\'{c}, Matschke, and Ziegler proved in \cite{blago-optimal-2011} the corresponding result for the Tverberg-Vre\'cica problem.  This effectively generalizes Theorem \ref{theorem-optimalcolored}.  The constraint method also yields results for the Tverberg-Vre\'cica conjecture \cite{blagojevic2016topological}.

Theorem \ref{theorem-dolnikovtransversal} is not optimal for a single set, and the {depth} of the affine subspace can be improved, as shown by Magazinov and P\'or \cite{magpor16} for $k=1$.  There is a nice conjecture by Bukh, Matou\v{s}ek, and Nivasch in this direction \cite{Bukh:2010hz}.

\begin{conjecture}\label{conjecture-boristransversal}
Let $0 \le k < d$ be integers.  Given a finite set $X$ of points in $\rr^d$, there is a $k$-dimensional affine subspace $L$ such that any closed half-space containing $L$ also has at least
\[
\frac{|X|(k+1)}{d+k+1}
\]
points of $X$.
\end{conjecture}

Conjecture \ref{conjecture-boristransversal} is known for $k\in \{0,d-2,d-1\}$.  This conjecture and the results by Magazinov and P\'or beg the question of whether the Tverberg-Vre\'{c}ica conjecture can be improved in the same way for a single set.  In other words, we present the following new conjecture.

\begin{conjecture}
Let $r,k,d$ be integers such that $0 \le k < d$.  Then, for any finite set $X$ of points in $\rr^d$ such that
\[
\left\lceil \frac{|X|(k+1)}{d+k+1}\right\rceil \ge r
\]
there is a partition of $X$ into $r$ sets $X_1, \ldots, X_r$ and a $k$-dimensional affine subspace $L$ such that $L$ intersects each of $\conv X_1 , \ldots, \conv X_r$.
\end{conjecture}

The case $k=0$ is Tverberg's theorem and the case $k=d-1$ follows easily by taking a halving hyperplane of $X$ for $L$ and pairing points of opposite sides of $L$ to form the partition. A {\sl halving hyperplane} is a hyperplane that has at least $|X|/2$ points of $X$ on both sides.  Note that if $|X|$ is odd, the halving hyperplane contains at least one of the points, which can be taken as a singleton in the partition.

Next we prove a version of the above conjecture with $k=d-2$ using $|X| = \left(\frac{2d-1}{d-1}\right)r+ O(d)$ points and the original method by Birch~\cite{Birch:1959}.  The case $k=d-2$ of Conjecture \ref{conjecture-boristransversal}, proved in \cite{Bukh:2010hz}, gives us an affine flat $L$ of dimension $d-2$ such that any half-space that contains it also has $\frac{|X|(d-1)}{2d-1} - O(d) = r$ points of $X$.

Notice that $L^{\perp}$ is a $2$-dimensional space and $L \cap L^{\perp}$ is a single point $p$.  Denote by $X^*$ the projection of $X$ onto $L^{\perp}$. We can order the points of $X^*$ clockwise around $p$ and  assign to them labels from $\{1,2,\ldots, r\}$ in such a way that if $x^*$ has label $i$, then the next point has label $i+1$ modulo $r$.  This gives us the partition of $X$ that we wanted.  Indeed, if all the points with label $j$ are separated strictly from $p$, then there is a closed half-plane $H^+$ that contains all of them but not $p$.  Then, the complement $H^-$ is an open half-plane that contains $p$ but has at most $r-1$ points of $X^*$, a contradiction.

Given a finite set of points in $\rr^d$, finding affine transversals to the convex hulls of its subsets in general is an interesting problem.  Consider the following instance.  Instead of seeking partitions that have a low-dimensional transversal, what if we seek a transversal to all sets of a given size?  Given $d,\lambda, k$ the following two parameters were introduced in \cite{Arocha:2011wo}.

First, $m(d,\lambda, k)$ is the maximum positive integer $n$ such that for any subset of $n$ points in $\rr^d$ there is an affine subspace of dimension $\lambda$ that intersects all the convex hulls of its subsets of cardinality $k$.

Second, $M(d,\lambda, k)$ is the minimum positive integer $n$ such that for every subset of $n$ points in general position in $\rr^d$ there is no affine subspace of dimension $\lambda$ that intersects all the convex hulls of the subsets of cardinality $k$.  The value of $M(d, \lambda, k)$ is known to be $(d-\lambda) + 2k +1 - \min\{k,\lambda\}$ \cite{Arocha:2011wo}, but the value of $m(d,\lambda, k)$ is still open.  The conjecture from Arocha et. al is the following (see also \cite{CMM17} for related results).

\begin{conjecture}
For $k,d$ positive integers and $0 \le \lambda \le d$ we have $m(d,\lambda, k) = (d-\lambda) + k + \lceil \frac{k}{\lambda}\rceil - 1$.
\end{conjecture}

\subsection{Reay's conjecture}

Tverberg's theorem gives us a partition of a set of $(r-1)(d+1)+1$ points into $r$ sets whose convex hulls all intersect.  If we instead ask for the convex hulls of every $k$ parts to intersect, it is not clear if a smaller number of points would be sufficient.  It was conjectured by Reay that this is not the case, even for $k=2$ \cite{Reay79}.

\begin{conjecture}
	There is a set $(r+1)(d-1)$ points in $\rr^d$ such that for any partition of them into $r$ parts, there are two parts whose convex hulls are disjoint.
\end{conjecture}

For general $k$ it brings the following problem.

\begin{problem}
Given positive integers $r,k,d$ such that $r \ge k \ge 2$, find the smallest integer $R(d,r,k)$ such that the following holds.  For any $R(d,r,k)$ points in $\rr^d$ there is a partition of them into $r$ parts $X_1, \ldots, X_r$ such that the convex hull of every $k$ of them intersect.	
\end{problem}

Reay's conjecture can be written as $R(d,r,k) = R(d,r,r)=(r-1)(d+1)+1$ for $k\ge 2$.  The best current general bound is $R(d,r,k)$ is $R(d,r,k) \ge r\left(\frac{k-1}{k} \cdot d + 1 \right)$ \cite{asada2016reay}.  Reay's conjecture is known to be true for $k \ge \frac{d+3}{2}$ or if $d < \frac{rk}{r-k}-1$, along a few other specific instances \cite{PeS16}.

\subsection{Tverberg with tolerance}

Tverberg's theorem also admits very robust versions, which resist the removal of points.  The first extension of this kind was proven by Larman \cite{Larman:1972tn}, also known as Radon's theorem with tolerance.

\begin{theorem}
	Given $2d+3$ points in $\rr^d$, there is a partition of them into two parts $A,B$ such that for any point $x$ we have
	\[
	\conv (A \setminus \{x\}) \cap \conv (B \setminus \{x\}) \neq \emptyset.
	\]
\end{theorem}

In other words, removing any single point won't break the Radon partition.  This result has been shown to be optimal for $1 \le d \le 3$ by Larman and for $d=4$ by Forge, Las Vergnas and Schuchert \cite{Forge:2001te}.  The best lower bound for this result is that at least $\left\lceil \frac{5d}{3} \right\rceil +3$ points are needed, which was proven by Ram\'irez-Alfons\'in using Lawrence oriented matroids \cite{Ram01}.

Extending Larman's result to partitions into more parts leads to a the following problem.

\begin{problem}\label{problem-tverbergtolerance}  Let $r,t,d$ be positive integers.  Determine the smallest integer $N=N(r,t,d)$ such that any set $X$ of $N$ points in $\rr^d$ can be partitioned into $r$ parts $X_1, \ldots, X_r$ such that for any set $C$ of at most $t$ points of $X$ we have
\[
\bigcap_{j=1}^r \conv (X_j \setminus C) \neq \emptyset.
\]	
\end{problem}

A surprising fact about this problem is that for fixed, $r,d$ we have $N = rt + o(t)$, which was first discovered by Garc\'ia-Col\'in, Raggi and Rold\'an-Pensado \cite{GRR17} using geometric Ramsey-type results.

The current best bounds for this result are as follows. $N = rt + \tilde{O}(\sqrt{t})$ for large $t$ and fixed $r,t$ \cite{soberon2016robust}, where the $\tilde{O}$ hides only polylogarithmic factors.  This bound is polynomial in all variables if the $\tilde{O}$ term is expanded, and is proved using the probabilistic method combined with Sarkaria's technique.

For small $t$, the bound above falls short of an earlier bound $N \le (r-1)(d+1)(t+1)+1$ \cite{Soberon:2012er}.  The only case when the optimal number is known is $d=1$, $N(r,t,1) = rt + 2r-1$, by Mulzer and Stein \cite{MS14}, who studied algorithmic versions of the problem.  Mulzer and Stein's bound for $d=2$, $N(r,t,2) \le 2(rt + 2r -1)$ is also the best known for some values of $r,t$.  For general lower bounds, using points in the moment curve gives $N(r,t,d) \ge r(t + \lfloor d/2 \rfloor +1 )$ \cite{soberon2015equal}.

The topological version of the problems with tolerance remains open, even the cases with $t=1$.  Tverberg with tolerance also has colorful versions, as in Conjecture \ref{conjecture-colored92}.  If we impose conditions on the partitions based on the colors, it is natural to also impose conditions on the points removed.  We define $N_{\operatorname{col}}(r,t,d)$ as the smallest integer such that, for any $N_{\operatorname{col}}$ sets of $r$ points each (considered as color classes), there is a partition of them into $r$ colorful sets $X_1, \ldots, X_r$ with the following property.  Even if we remove any $t$ color classes, the convex hulls of what is left in each $X_j$ still intersect.  It is known that for $r,d$ fixed, and $r \ge 3$ we have $N_{\operatorname{col}}(r,t,d) \le t(1.6 +o(1))$, and $N_{\operatorname{col}}(2,t,d) \le t(2+o(1))$ \cite{soberon2016robust}.   However, it may be that fewer color classes are needed.

  \begin{conjecture}
  For $r,d$ fixed, we have $N_{\operatorname{col}}(r,t,d) = t(1+o(1))$.
  \end{conjecture}
  
If we want to remove a larger proportion of points while still having a Tverberg partition, we need several partitions.  The number of partitions needed was determined in \cite{soberon2017epsilon}.

\begin{theorem}
Let $\varepsilon >0$ be a real number and $r,d,m$ be positive integers such that $\varepsilon > (1-1/r)^m$.  Then, for all finite sets $X \subset \rr^d$ of sufficiently large cardinality, we can find $m$ partitions of $X$ into $r$ sets each, such that for any subset of at least $\varepsilon |X|$ points of $X$, at least one of the partitions induces a Tverberg partition.
\end{theorem}

The condition on $\varepsilon$ is sharp.  This result follows from extending the probabilistic approach of \cite{soberon2016robust}.  The case $m =1$ is essentially the statement $N(r,t,d) = rt + o(t)$.  It can be interpreted as a version with tolerance $(1-\varepsilon)|X|$.

\subsection{Dimension of Tverberg points}

Assume $X \subset \rr^d$, and $r$ is a positive integer, and define $T_r(X)$ as \textit{the set of points which are in the intersection of the convex hulls of some $r$-Tverberg partition of $X$.} Thus $T_r(X)$ is simply the union of all $\bigcap_1^r \conv X_i$ taken over all $r$-partition of $X$. The size or dimension of $T_r(X)$ is another way to quantify Tverberg partitions. Here we do not assume that $X$ is in general position, $T_r(X)$ is interesting even in that case. 

For a set $A \subset \rr^d$, we consider $\dim (A)$ the Hausdorff dimension of $A$, with the convention $\dim(\emptyset) = -1$.  With this definition the following conjecture was made by Kalai in 1974 \cite{Kal00flavor}.

\begin{conjecture}[Cascade conjecture]\label{conjecture-kalai}
For any finite set $X$ of points in general position in $\rr^d$ we have
\[
\sum_{s=1}^{|X|} \dim(T_s(X)) \ge 0.
\]
\end{conjecture}

\begin{figure}
\centerline{\includegraphics[scale=0.8]{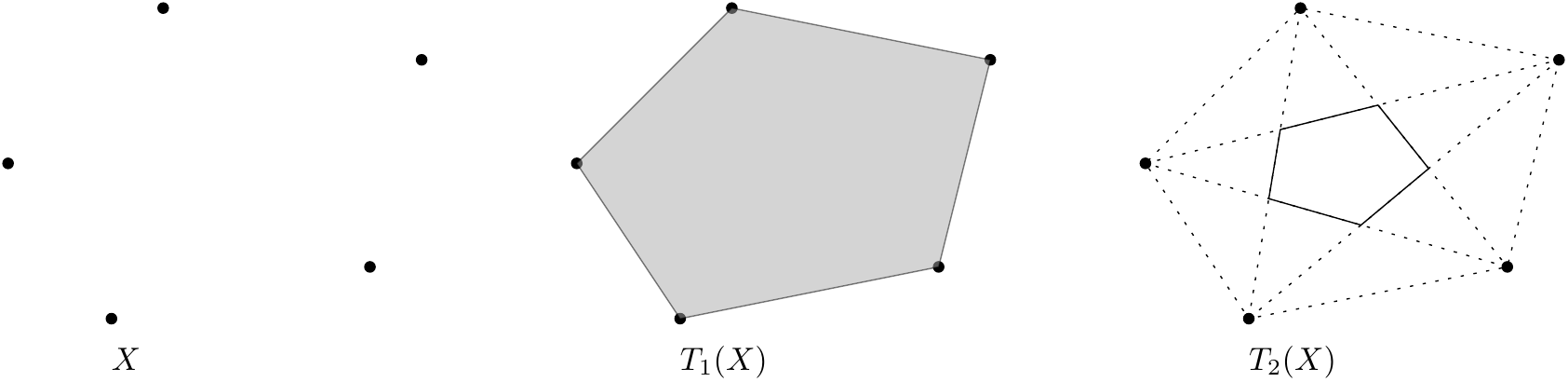}}
\caption{This figure shows $\dim(T_1(X))=2$ and $\dim(T_2(X))=1$ for a set $X$ of five points in the plane.  For this set, $T_3(X), T_4(X), T_5(X)$ are empty.  Notice that $\sum_{s=1}^5 \dim(T_s(X)) = 2 + 1 + (-1) + (-1) + (-1) = 0$, as expected from Conjecture \ref{conjecture-kalai}.  The example also agrees with Conjecture \ref{conjecture-reaydimension} with $d=2, r=2, k=1$.}
\end{figure}

This conjecture implies Tverberg theorem. Indeed, if $|X| = (r-1)(d+1)+1$ and has no $r$-Tverberg partitions, then $\dim(T_s(X)) \le d$ for $s < r$ an $\dim (T_s(X)) = -1$ for all $s \ge r$, which means that the sum above would be at most $-1$.  A weaker version is also open.

\begin{conjecture}[Weak cascade conjecture]
For any finite set $X$ of points in general position in $\rr^d$ we have
\[
\sum_{s=1}^{|X|} \dim(\conv T_s(X)) \ge 0.
\]
\end{conjecture}
Here $\dim(\conv T)$ is the usual dimension of $\conv T$ or $\aff T$. Another conjecture in this direction was formulated by Reay \cite{Reay68} shortly after Tverberg published his result.

\begin{conjecture}\label{conjecture-reaydimension}
Let $k,d,r$ be integers such that $0 \le k \le d$.  Then, for any set $X$ of $(r-1)(d+1)+k+1$ points in $\rr^d$ in general position we have
\[
\dim(T_r(X)) \ge k.
\]
\end{conjecture}

The case $k=0$ is Tverberg's theorem.  Reay proved his conjecture when the points are in \textit{strongly} general position, but believed that just general position (i.e. no $d+1$ points of $X$ lie in a hyperplane) should be enough.  The conjecture has been proved for $d \le 8$ and any $k,r$, for $r \le 8$ and any $d,k$ and for $d=24, k=1$ and any $d$, see \cite{Roudneff:2001cl, Roudneff2001part2, Roudneff09}.

We note that for sets of points in general position, Conjecture \ref{conjecture-kalai} and Conjecture \ref{conjecture-reaydimension} are equivalent.

{\bf Proof.}
Assume first that Conjecture \ref{conjecture-reaydimension} holds, and let $X$ be a set of $n$ points in general position in $\rr^d$.  Then, there are non-negative integers $r,k$ such that $n = (d+1)(r-1)+k+1$ with $k \le d$.  Then, by Conjecture \ref{conjecture-reaydimension} we have $\dim((T_s(X)) = d$ for $s \le r-1$, $\dim(T_r(X)) \ge k$ and $\dim((T_s(X)) = -1$ for $s \ge r+1$ by the general position assumption.  Therefore

\begin{align*}
\sum_{s=1}^{|X|} \dim( T_s(X)) & \ge d(r-1) + k + (-1)(n-r) \\
& = (d+1)(r-1)+ k + 1 - n = 0.
\end{align*}

If we now assume \ref{conjecture-kalai}, and $X$ is a set of $(d+1)(r-1)+k + 1$ points in general position, note that $\dim(T_s(X)) \le d$ for $s < r$, and $\dim(T_s(X)) = -1$ for $s> r$ by the general condition assumption.  Therefore

\begin{align*}
0 & \le \sum_{s=1}^{|X|} \dim( T_s(X)) \\
& \le d(r-1) + \dim(T_r(X)) + (-1)(|X|-r) \\
& = \dim(T_r(X)) - k. 
\end{align*}
\qed

This implies that Kalai's cascade conjecture holds for sets of points in sufficiently general position, and for points in general position where Roudneff has proven Reay's conjecture.

\subsection{Finding Tverberg partitions}


Centerpoints and notions of depth are key concepts in data analysis.  Centerpoints often play the role of high-dimensional median.  For an $n$-point set in $\rr^d$, there are algorithms that find a centerpoint in time $O(n^{d-1})$, which is believed to be optimal \cite{Chan04center}.

In general, it is computationally difficult to verify the depth of a point in data set.  However, given a Tverberg partition $X_1, \ldots, X_r$ with a point $p$ in the convex hull of each $X_j$,the depth of $p$ is at least $r$, clearly. This a lower bound can be verified in polynomial time: simply check that $p \in \conv X_j$ for each $j$.  This implies that $r$ is a lower bound for the depth of $p$ in the set of points.

Given a set of $n$ points in $\mathbb{R}^d$, finding Tverberg partitions into $\left\lceil \frac{n}{d+1}\right\rceil$ parts in polynomial time is out of reach for current algorithms, and an interesting open problem by itself.  To achieve fast algorithms, we have to pay the price of reducing the number of parts in our partition.  There is a deterministic algorithm by Miller and Sheehy that gives a Tverberg partition with $r = \left\lceil \frac{n}{(d+1)^2}\right\rceil$ in $n^{O(\log d)}$ time \cite{MSR10approx}.  Using a lifting argument in combination with Miller and Sheehy's algorithm, a deterministic algorithm that gives a Tverberg partition with $r = \left\lceil \frac{n}{4(d+1)^3}\right\rceil$ that runs in time $d^{O(\log d)}n$ was produced by Mulzer and Werner \cite{MD13linear}.  This is linear in $n$ for any fixed dimension.

For non-deterministic arguments, one can compute Tverberg points with $r = \left\lceil \frac{n}{d(d+1)^2}\right\rceil$ with a probability $\varepsilon>0$ of failure fixed in advance \cite{RS16center}.  This algorithm is weakly polynomial in all variables $n,d,\log(1/\varepsilon)$.

The algorithmic versions of other variations of Tverberg's theorem are also interesting.  For instance, for Tverberg's theorem with Tolerance, Mulzer and Stein showed how the two deterministic algorithms described above could be adapted to that setting \cite{MS14}.  If one is willing to have non-deterministic arguments, the results in \cite{soberon2016robust} show that by randomly assigning each point to one of $X_1, \ldots, X_r$ independently, we can bound the probability of failure efficiently. 

\section{Universal Tverberg partitions}\label{section-universal}

Assume $X \subset \R^d$, and $|X|=(r-1)(d+1)+1$. A natural question is which $r$-partitions of $X$ are Tverberg partitions. One case when this structure is completely known is when the points of $X$  come from the moment curve $m(t)=(t,t^2,\ldots,t^d) \in \R^d$, $t \ge 0$ and are far apart from each other. Note that the points $m(t_1),m(t_2),\ldots,m(t_N)$ on the moment curve are ordered by $t_1<t_2<\ldots<t_N$. In order to understand the structure of Tverberg partitions it is better to work with sequences $(a_1,\ldots,a_N)$ of points in $\rr^d$ (instead of sets $X\subset \rr^d$). So in this section we work with sequences. 

We start with the simplest case: that of Radon partitions. Consider $d+2$ points $m(t_1),\ldots,m(t_{d+2})$ on the moment curve with $t_1<\ldots<t_{d+2}$. It is well-known (see for instance Gr\"unbaum's book~\cite{Grun} or Matou\v{s}ek's~\cite{Matousek:2002td}) that there is a unique Radon partition in this case, namely, one set is $X_1=\{m(t_i): i \mbox{ odd}\}$ and the other one is $X_2=\{m(t_i): i \mbox{ even}\}$. That is, the Radon partition is just two {\sl interlacing} sets, meaning that on the moment curve between two consecutive points of $X_1$ (resp $X_2$) there is a point of $X_2$ (and $X_1$). It is also known that this is the {\sl universal} Radon partition: for every $d \in \N$ there is $N\in \N$ such that any $d$-dimensional (general position) vector sequence $a_1,\ldots,a_N$ contains a subsequence $a_{i_1}\ldots,a_{i_{d+2}}$ with $i_1<i_2<\ldots<i_{d+2}$ such that their unique Radon partition is the interlacing sets $X_1=\{a_{i_j}: j \mbox{ odd}\}$ and $X_2=\{a_{i_j}: j \mbox{ even}\}$. The moment curve shows that this is the unique universal Radon partition. 

\smallskip
What is the corresponding statement for Tverberg partitions? 

\begin{figure}
\centerline{\includegraphics[scale=1]{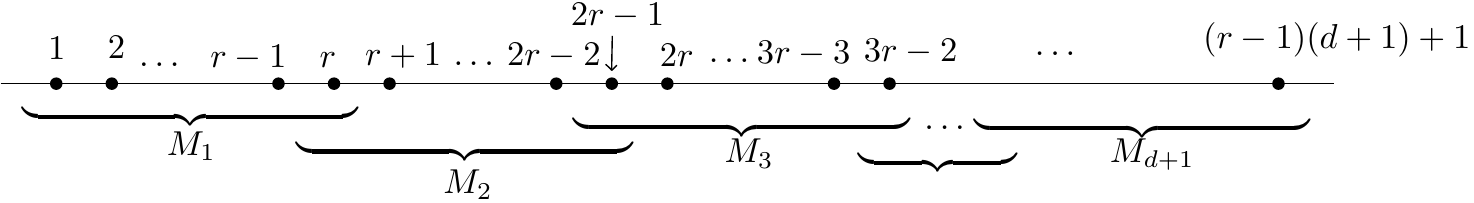}}
\caption{All blocks $M_i$ are of length $r$ and share one element with $M_{i-1}$. In a special $r$-partition of $[n]$ $(n= (r-1)(d+1)+1)$, all the $r$ elements of each block are in different parts of the partition.}
\label{figure-blocks}
\end{figure}

\smallskip
We need some definitions for partitions of sequences. For $k\in [d+1]$ we define the {\sl block} $M_k$ as the set of consecutive integers $M_k=\{(r-1)(k-1)+1,(r-1)(k-1)+2,\ldots,rk-1\}$. The blocks almost form a partition of $[n]$ with $n=(r-1)(d+1)+1$, only the elements $r,2r-1,3r-2,\ldots,rd-(d-1)$ are covered twice, namely by $M_1,M_2$, $M_2,M_3$, etc, $M_d,M_{d+1}$.  See Figure \ref{figure-blocks}. We Call an $r$-partition $I_1,\ldots,I_r$ of $[n]$ {\sl special} if $|I_j\cap M_k|=1$ for every $j\in [r]$ and every $k\in [d+1]$. 

\begin{figure}
\centerline{\includegraphics[scale=1]{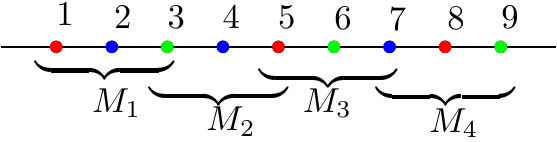}}
\caption{A special partition for $d=3, r=3$.  In this example $I_1 = \{1,5,8\}, I_2=\{2,4,7\}, I_3 = \{3,6,9\}$.  Note that each $I_j$ has exactly one element of each $M_i$.}
\end{figure}

Now let $0<t_1<\ldots <t_n$ be a  {\sl rapidly increasing} sequence of real numbers, meaning that, for every $h \in [n-1]$, $t_{h+1}/t_h$ is at least as large as some (large) constant $c_{d,r,h}$ depending only on $d,r,h$.  Consider the set $X$ of points in $\rr^d$ in the moment curve, $X = \{m(t_1), \ldots, m(t_n)\}$.

It is clear that a partition $P=\{I_1,\ldots,I_r\}$ of $[n]$ induces the partition $Q=\{X_1,\ldots,X_r\}$ of $X$ via $X_j=\{m(t_i): i\in I_j\}$ and vice versa. Note that when $r=2$, that is Radon partitions, each block contains two elements and the special partitions are exactly the interlacing ones. It was observed by B\'ar\'any and P\'or, and by Mabillard and Wagner (both unpublished), that the Tverberg partitions into $r$ parts of the point set $X=\{m(t_1),\ldots,m(t_n)\}$ can be explicitly described: assuming that the constants $c_{d,r,h}$ are suitably large,  $Q$ is a Tverberg partition of $X$ if and only if $Q$ is induced by a special partition $P$ of $[n]$.

It is known (and easy to check) that the number of special partitions of $[n]$, which is then the same as the number of Tverberg partitions of $X$, is equal to $(r-1)!^d$ when $X=\{m(t_1),\ldots,m(t_n)\}$ and the sequence $t_1,\ldots,t_n$ is rapidly increasing. So this is another example achieving the lower bound in Skiersma's conjecture.

A similar example was given by Bukh, Loh, and Nivasch in~\cite{bukh2016classifying}, but instead of the moment curve, they use the ``diagonal of the stretched grid'',  for the exact definition see ~\cite{bukh2016classifying}. The example is again a sequence $X=\{a_1,\ldots,a_n\}$ of points in $\R^d$ with the property that the partition $Q=\{X_1,\ldots,X_r\}$ of $X$ is a Tverberg partition if and only if the corresponding partition $P=\{I_1,\ldots,I_r\}$ of $[n]$ is special. Here, of course, $X_j=\{a_i:\;i\in I_j\}$ again. The example also achieves the lower bound in Sierksma's conjecture, as one can easily check.

Which Tverberg partitions must always appear in every set $X\subset \R^d$ of large enough cardinality (when $X$ is in general position)?  In other words, consider an $r$-partition $P$ of $[n]$.  Given a sequence $a_1, \ldots, a_N$ of $N$ points in $\rr^d$ with $N$ large enough, can one always find a subsequence $b_1, b_2,\ldots, b_n$ where $P$ induces a Tverberg partition? If such a partition always exists, then $P$ is called {\sl unavoidable}: the name indicates that such a partition is always present in a long enough sequence. The previous examples show that an unavoidable partition has to be a special partition of $[n]$. Bukh, Loh and Nivasch showed in~\cite{bukh2016classifying} that in low dimensions  ($d=1,2$) the unavoidable partitions are exactly the special ones, and conjectured that in high dimensions the same statement holds. Shortly after, P\'or characterized such partitions $P$ in all dimensions \cite{pornew}, answering positively the conjecture in \cite{bukh2016classifying}. Actually, P\'or proved the following {\sl universality theorem} for Tverberg partitions of vector sequences.

\begin{theorem}[Universal Tverberg partitions]\label{theorem-por} Given integers $d\ge 1, r\ge 2$ and $m\ge n=(r-1)(d+1)+1$, there is an integer $N$ with the following property. Every sequence $a_1,\ldots,a_N \in \R^d$ of vectors in general position contains a subsequence $b_1,\ldots,b_m$ such that for every subsequence $b_{i_1},\ldots,b_{i_n}$ with $1\le i_1<i_2 < \ldots < i_n\le m$, the Tverberg partitions are exactly the ones induced by the special partitions of $[n]$.
\end{theorem}

\section{Applications of Tverberg's theorem}\label{section-applications}

An early application, actually the motivation for both Birch \cite{Birch:1959} and Tverberg, was Rado's centerpoint theorem, described in Section 1.1. The set $C(X)$ of center points for a given set $X \in \R^d$ is a convex set that contains every Tverberg point. It is known, however, that the convex hull of the Tverberg points does not coincide with $C(X)$  \cite{Avis:1993jk}.

Another geometric application is the so called {\sl first selection lemma}. It states the following.

\begin{theorem}\label{th:select} Given a set $X$ of $n$ points in $\R^d$ (in general position), there is a point $z\in \R^d$ which is contained in the convex hull of at least $c_d{n \choose {d+1}}$ of the ${n \choose {d+1}}$ possible $(d+1)$-tuples of $X$, here $c_d\ge (d+1)^{-d}$ is a constant depending only on $d$.
\end{theorem}

This result was proved by Boros and F\"uredi \cite{Boros:1984ba} for $d=2$, the general case by B\'ar\'any~\cite{Barany:1982va}. The exact value of the constant $c_d$ is known only for $d=2$ \cite{Bukh:2011vs} and  \cite{Boros:1984ba}, where $c_2 = 2/9$. The proof in \cite{Barany:1982va} shows $c_d\ge (d+1)^{-d}$. It is known \cite{Bukh:2011vs} that $c_d \le \frac {(d+1)!}{(d+1)^{d+1}}$. There was a slight improvement by Wagner~\cite{Wagner:2003us}. In a remarkable paper Gromov~\cite{Gromov:2010eb} gives an exponential improvement by showing that
\[
c_d \ge \frac {2d}{(d+1)(d+1)!}\sim \frac {e^d}{(d+1)^{d+1}}.
\]
In fact, Gromov proves the following stronger, topological statement: for every continuous map $f:\skel_d \tri^{n-1} \to \rr^d$ there is a point in $\R^d$ whose preimage intersects at least
\[
\frac {2d}{(d+1)(d+1)!}{n \choose {d+1}}
\]
faces of dimension $d$. Theorem~\ref{th:select} is the special case when $f:\skel_d \tri^{n-1}\to\rr^d$ is an affine map. Surprisingly, the topological proof gives a better constant.  A simplified proof appeared in \cite{Karasev:2012bj}.

\begin{problem} What is the order of magnitude of the constant $c_d$?  Does $(d+1)^{d+1}c_d$ exhibit exponential or superexponential growth? Also, are the constants for the topological and affine versions of the problem equal?
\end{problem}


One seminal application of Tverberg's theorem is the weak $\varepsilon$-net theorem for convex sets \cite{Alon:1992ek}.

\begin{theorem}\label{theorem-nets}
Let $d$ be a positive integer and $\varepsilon >0$ a real number.  Then, there is a constant $n=n(d,\varepsilon)$ such that for each finite set $X$ of points in $\rr^d$, there is a set $P$ of $n(d,\varepsilon)$ points such that for each $Y \subset X$ with $|Y|\ge \varepsilon |X|$, we have
\[
P \cap \conv Y \neq \emptyset.
\]
\end{theorem}

This is a very strong result on the combinatorial properties of convex sets.  The reader may verify that the equation $n(d, d/(d+1)) = 1$ is the centerpoint theorem.  The weak $\varepsilon$-net theorem for convex sets is proved by repeatedly using the first selection lemma to greedily construct the set $P$, one point at a time.  If one applies Gromov's topological extension of Theorem \ref{th:select} instead of the first selection lemma, we obtain a topological version of Theorem \ref{theorem-nets} \cite{MS17positive}.  The weak $\varepsilon$-net theorem was a key component of the proof of Hadwiger-Debrunner $(p,q)$ conjecture \cite{Alon:1992ta} (cf. \cite{Alon:1992ek}), a celebrated result in combinatorial geometry.

Another application of Tverberg's theorem, or rather of its colorful version, concerns halving planes. In this case Tverberg's theorem helped to locate the key question in the following way. A halving plane of a finite set $X \subset \R^3$ of points in general position is a plane spanned by three points of~$X$ that has equally many points of $X$ on either side of it. (So $|X|=n$ is odd.) While the first author was working with F\"uredi and Lov\'asz on establishing upper bounds on the number of halving planes they encountered the following question: given a set $X \subset \R^2$ of $n$ points in general position, a {\sl crossing} is the intersection of the lines spanned by $x,y$ and by $u,v$ where $x, y, u, v$ are distinct points from~$X$. It is evident that there are $\frac12\binom{n}{2}\binom{n-2}{2} \sim n^4$
crossings. How many of them are contained in a typical triangle spanned by points in~$X$?
A direct application of Tverberg's theorem combined with a double counting argument shows that the number of crossings is again of order~$n^4$. This was the first step in establishing an $O(n^{3-\varepsilon})$ bound on the number of halving planes. The proof uses the supersaturated hypergraph lemma of Erd\H os and Simonovits~\cite{ErdSim}, and that is why a special version of Tverberg's theorem, a colorful variant was needed cf.~\cite{Barany:1990wa}. The method was extended to higher dimensions in \cite{Alon:1992ek}.  This is how the halving plane question lead to the colorful Tverberg theorem. The moral is that when working on a question in combinatorial convexity it is always good to check what Tverberg's theorem says in the given situation.

Another application of the colorful Tverberg theorem is a result by Pach~\cite{Pach:1998vx} on homogeneous selection:

\begin{theorem}\label{th:pach} Assume we are given sets $C_1,\ldots,C_{d+1}\subset \R^d$ (considered as colors classes) that have the same size $|C_i|=n$ for all $i \in [n]$. Then there are subsets $Q_i \subset C_i$ with $|Q_i|\ge c_d n$ and a point $z \in \R^d$ such that $z \in \conv\{x_1,\ldots,x_{d+1}\}$ for every transversal $x_1 \in Q_1,\ldots,x_{d+1} \in Q_{d+1}$. Here $c_d>0$ is a constant depending only on $d$.
\end{theorem}

There are further geometric applications of Tverberg's theorem in \cite{bukh2017one} and in quantum correcting codes~\cite{MR1745959}.

Here is a purely combinatorial result, originally a theorem of Lindstr\"om~\cite{Lin72} that turned out to be a consequence of Tverberg's theorem.

\begin{theorem}\label{th:lindst} Assume $n,r>1$ are integers and set $N=(r-1)n+1$. If $A_1, \dots, A_N$ are non-empty subsets of an $n$-element set, then there are non-empty and disjoint subsets $J_1, \dots, J_r$ of $[N]$ such that $\bigcup_{i \in J_1} A_i = \dots = \bigcup_{i \in J_r} A_i$.
\end{theorem}

The geometric proof, found by Tverberg himself \cite{tverbergLind}, transfers this purely combinatorial partition problem to convex geometry.

{\bf Proof.} We assume that the ground set is $[n]$, that is $A_i \subset [n]$. Associate with each set $A_i$ the vector
\[
a_i=\frac {\chi_{A_i}}{|A_i|}
\]
where $\chi_{A_i}$ is the characteristic vector of $A_i$. So $a_i$ is in $\R^n$ but it lies, in fact, in the affine subspace $S$ where the sum of the coordinates is equal to one. This subspace is a copy of $\R^{n-1}$.  We can apply Tverberg's theorem to the points $a_1,\ldots,a_N\in S$.  This gives us a partition $I_1,\ldots,I_r$ of $[N]$ and a point $a \in S$ with
\[
a \in \bigcap_{h=1}^r \conv \{a_i: i \in I_h\}.
\]
The common point $a\in S \subset \R^n$ is a non-zero vector in $\R^n$ with non-negative coordinates. Let $J\subset [n]$ be the set of non-zero coordinates of $a$. It is easy to see that for a suitable subset $J_h$ of $I_h$, $J$ is the union of $A_i$ with $i \in J_h$ for every $h\in [r]$.
\qed

The result above can be extended to bound the number of such partitions, effectively proving the analogue of Sierksma's conjecture in tropical geometry \cite{tropicaltverberg}.

There is a recent and very powerful application of the topological Tverberg theorem, due to Frick~\cite{FFrick1,  FFrick2}.  One of the first examples of a combinatorial problem solved with topological methods is Lov\'asz's groundbreaking proof~\cite{Lovasz} of Kneser's conjecture.  He used the Borsuk-Ulam theorem to establish a lower bound on the chromatic number of Kneser graphs. Since then, topological methods have been used to bound the chromatic number of graph and hypergraphs. 

The connection between Tverberg type results and Kneser hypergraphs was first noted and used by Sarkaria in 1990 \cite{Sarkaria:1990} and \cite{Sarkaria:1991ug}. Going much further, Frick elucidates the underlying connection between intersection patterns of finite sets and topological statements, via Tverberg-type theorems. This creates a dictionary between the two types of results.  To state just one theorem of this kind, let $L$ be a simplicial complex and $K \subset L$ be a subcomplex. Denote by $KN^r(K,L)$ the $r$-uniform hypergraph whose vertices are the inclusion-minimal faces of $L$ that are not contained in $K$ and whose hyperedges are the $r$-tuples of vertices when the corresponding faces are pairwise disjoint.

\begin{figure}
\centerline{\includegraphics[scale=0.8]{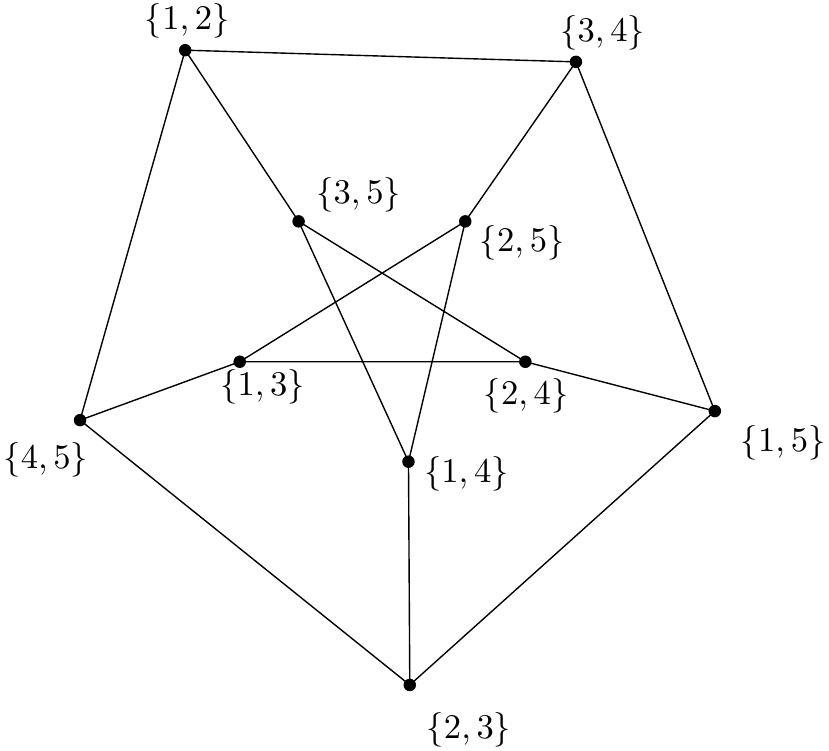}}
\caption{A labeled figure of $KN^2(\skel_0{\Delta^4}, \Delta^4)$.  The vertices are pairs of integers in $[5]$ and there is an edge between two pairs if they are disjoint.}
\label{figure-petersen}
\end{figure}

For example, if $r=2$, $L= \Delta^4$ and $K=\skel_0 \Delta^4 = [5]$, then $KN^r(K,L)$ is the 1-skeleton of $\Delta^4$, that is, the Petersen graph (see Figure \ref{figure-petersen}).  If $r=2$, $L=\tri^{n-1}$ and $K = \skel_{k-2} \Delta^{n-1}$, then the vertices of $KN^2(K,L)$ are the $k$-tuples of $[n]$, with two connected if they are disjoint. This is exactly  the Kneser graph of $k$-subsets of $[n]$.

The general principle behind the constraint method developed in \cite{Blagojevic:2014js} is then used to prove the following result~\cite{FFrick1}.

\begin{theorem}\label{th:frick}
Assume $d,k\ge 0$ and $r\ge 2$ are integers.  Let $L$ be a simplicial complex such that for every continuous map $g: L \to \R^{d+k}$ there exist disjoint faces $\sigma_1,\ldots,\sigma_r$ of $L$ such that $g(\sigma_1)\cap\ldots\cap g(\sigma_r)\ne \emptyset$.  If $\chi(KN^r(K,L))\le k$ for some subcomplex $K$ of $L$, then for every continuous map $f:K \to \R^d$ there are $r$ pairwise disjoint faces $\sigma_1,\ldots,\sigma_r$ of $K$ such that $f(\sigma_1)\cap\ldots\cap f(\sigma_r)\ne \emptyset$.
\end{theorem}

This result can be used in two directions: establishing the upper bound $\chi(KN^r(K,L))\le k$ proves the existence of an $r$-fold intersection point for every continuous map $f:K \to \R^d$, and exhibiting a  continuous map $f:K \to \R^d$ without such an $r$-fold intersection gives the lower bound  $\chi(KN^r(K,L))\ge k+1$. The theorem relates intersection patterns of continuous images of faces in a simplicial complex to intersection patterns of finite sets. It implies, generalizes, and unifies several earlier results of this type, including those by Lov\'asz~\cite{Lovasz}, Dol'nikov~\cite{Doln}, Alon, Frankl, Lov\'asz~\cite{AlonFrLo} and K\v{r}\'i\v{z}~\cite{kriz}.

\section{Tverberg-type results in distinct settings}\label{section-spaces}

\subsection{Convexity spaces and $S$-convexity}
Many variations of Tverberg's theorem appear if we change the underlying space we are using.  For example, consider the following integer version of Tverberg's theorem.

\begin{theorem}[Integer Tverberg]
Given $r,d$ positive integers, there is an integer $L=L(r,d)$ such that for any set of $L$ points in $\rr^d$ with integer coordinates there is a partition of the set into $r$ parts $L_1, \ldots, L_r$ such that the intersection of their convex hulls contains a point with integer coordinates.
\end{theorem}

The exact values of $L(r,d)$ are not known, even for $r=2$.  The existence of $L(r,d)$ follows from results by Jamison \cite{Jamison:1981wz}, as was noted by Eckhoff \cite{Eckhoff:2000jw}.  The number of points needed is much larger than in the non-integer case.  For instance, we have the lower bound $L(r,d) > (r-1)2^d$.  To see this, take $r-1$ copies of each vertex of the hypercube $[0,1]^d$.  This lower bound may not always be optimal.  For $r=2$, Onn proved $\frac{5}{4}2^d + 1 \le L(2,d) \le d(2^d-1) +3$ \cite{Onn:1991er}.  The case $d=3$ remains interesting, with $L(2,3) \le 17$ being the best upper bound \cite{BB03}. The best general upper bound up to date is $L(r,d) \le (r-1) d 2^d + 1$ \cite{discrete-tverberg17}.

\begin{problem}
Determine the value $L(2,3)$, or improve the current bounds $11 \le L(2,3) \le 17$.	
\end{problem}

To state properly Tverberg's theorem in abstract terms we only need two ingredients.  First is $\C_d$, the family of all sets in $\rr^d$ that are considered convex, and the second is to be able to compute convex hull; i.e., an operator $\conv : 2^{\rr^d} \to \C_d$ with a few properties.  Thus, given a ground set $Y$, a way to axiomatize convexity is to have an operator $\conv: 2^Y \to 2^Y$ which satisfies the following.

\begin{itemize}
	\item $\conv ( \conv A)) = \conv A$ for all $A \subset Y$,
	\item $A \subset \conv A$ for all $A \subset Y$,
	\item $A \subset B \subset Y$ implies $\conv A \subset \conv B$,
	\item For a countable sequence $A_1 \subset A_2 \subset \ldots \subset Y$ we have that $\cup_{i=1}^{\infty} \conv (A_i) = \conv \left( \cup_{i=1}^{\infty} A_i\right)$.
\end{itemize}

We say that the pair $(Y, \conv)$ is a \textit{convexity space}.  A central questions of convexity spaces is the following.

\begin{problem}
Given a convexity space $(Y, \conv)$ and a positive integer $t$, determine the value of $r_{t}$ (if it exists) such that for any $X \subset Y$ of $r_t$ points there is a partition of $X$ into $t$ parts $X_1, \ldots, X_t$ such that
\[
\bigcap_{j=1}^t \conv X_j \neq \emptyset.
\]
\end{problem}

The reason for the conflicting notation with our use of the variable $r$ is that, in the context of convexity spaces, the number $r_t$ above is called \textit{the $t$-th Radon number}.  A classic conjecture by Eckhoff was that for any convexity space with a finite $r_2$ we have $r_t \le (t-1)(r_2-1)+1$.  In other words, it asked if Tverberg's theorem follows from Radon for purely combinatorial reasons. Hopes for this were dashed by an example, presented by Boris Bukh in an unpublished preprint, that constructs a convexity space with $r_2 = 4$ and $r_t \ge 3(t-1)+2$ \cite{bukh2010radon}.  It remains open whether $r_t$ can be bounded as a function that is linear in both $r_2$ and $t$, which is enough for several applications.

General convexity spaces are outside the scope of this survey.  The interested reader should consult Eckhoff's survey \cite{Eckhoff:2000jw} on the subject, which also discusses convexity spaces where Eckhoff's conjecture is known to hold.  We focus on convexity spaces which are closely related to convexity in $\rr^d$.

An example is $S$-convexity, which generalizes the integer case.  Given a set $S \subset \rr^d$, we say that a set $A \subset S$ is $S$-convex if $A = S \cap \conv A$, where $\conv(\cdot )$ denotes the usual convex hull in $\rr^d$.  Given $A \subset S$, we define the $S$-convex hull $\conv_S(A)$ as the intersection of all $S$-convex sets $B$ such that $A \subset B$.

A Tverberg-type theorem for $S$ would simply be the existence of a number $T_S (r,d)$ such that for any set $X$ of $T_S(r,d)$ points of $S$, there is a partition of them into $r$ sets $X_1, \ldots, X_r$ such that
\[
\bigcap_{j=1}^r \conv_S(X_j) \neq \emptyset.
\]
It turns out that the existence of such theorems relies on whether there is a Helly-type theorem for $S$-convexity \cite{discrete-tverberg17}.  An $S$-convexity Helly theorem simply says that \textit{there is a natural number $h(S)$ such that, for any finite family of convex sets in $\rr^d$, if the intersection of any $h(S)$ or fewer of them contains a point of $S$, then the intersection of the whole family contains a point of $S$.}  The existence of Tverberg-type theorems is given by the following theorem.

\begin{theorem}
If $\conv_S$ has a Helly theorem with Helly number $h(S)$, then it has a Tverberg theorem.  Moreover, $T_S(r,d) \le h(S) d (r-1) + 1$ for all $r$.
\end{theorem}

This result implies the upper bound for integer Tverberg if we use Doignon's theorem, which says that $h(\mathbb{Z}^d) = 2^d$ \cite{Doignon:1973ht, Scarf:1977va, Bell:1977tm}.  If we have a \textit{quantitative} Helly for $S$, i.e., a natural number $h_k(S)$ such that, \textit{for any finite family of convex sets in $\rr^d$, if the intersection of any $h_k(S)$ of them contains at least $k$ points of $S$, then the intersection of the whole family contains at least $k$ points of $S$}, then we also obtain a similar Tverberg theorem, where now $\bigcap_{j=1}^r\conv_S(X_j)$ contains at least $k$ points of $S$.  The related Helly-type results on $S$-convexity are described in \cite{Amenta:2015tp}.

Although $S$-convexity often gives worse bounds than the classic setting, it is interesting that in some cases the asymptotic behavior of variations of Tverberg's theorem remains the same.  For example, we can naturally ask for a version of Tverberg with tolerance, Problem \ref{problem-tverbergtolerance}, makes sense in the integer lattice.  We get the following result.

\begin{theorem}
	Let $r,t,d$ be positive integers, where $r,d$ are fixed.  Then, there is a number $L(t) = rt + o(t)$ such that the following holds.  Given a set $X$ of $L(t)$ points with integer coordinates in $\rr^d$, there is a partition of $X$ into $r$ sets $X_1, \ldots, X_r$ such that for any for any $C \subset X$ of cardinality at most $t$, the convex hulls $\conv( X_1 \setminus C), \ldots, \conv(X_r \setminus C)$ have a common point with integer coordinates.
\end{theorem}

For a proof, we only need to follow the methods of \cite{GRR17} verbatim.  When they require the usage of a centerpoint, we simply need to use an integer point of depth ${2^{-d}}$, which exist as a consequence of Doignon's theorem.  It is interesting that for this version the linear-algebraic methods that use Sarkaria's technique fail completely.

\subsection{Tverberg-type theorem on families of sets.}

It is possible to prove Tverberg-type results if we are dealing with families of subsets of $\rr^d$ instead of just points.  In this case, we have to replace the convex hulls by other operators, or require a conclusion stronger than the intersection of the convex hulls being non-empty.  For example given a family $\mathcal{F}$ of $d+1$ hyperplanes in general position, we denote by $\tri(\mathcal{F})$ the simplex whose faces are given by $\mathcal{F}$.  Then, we have the following Tverberg-type result by Karasev \cite{Kar08hyp, Karasev:2011jv}.

\begin{theorem}\label{theorem-karasevhyperplanes}
Let $r,d$ be positive integers such that $r$ is a prime power.  Let $\mathcal{F}$ be a family of $r(d+1)$ hyperplanes in general position in $\rr^d$.  Then, there is a partition of $\mathcal{F}$ into $r$ sets $\mathcal{F}_1, \ldots, \mathcal{F}_r$ of $d+1$ hyperplanes each such that
\[
\bigcap_{j=1}^r \tri(\mathcal{F}_j) \neq \emptyset
\]
\end{theorem}

For this result there is also a corresponding discrete version of a centerpoint theorem for hyperplanes.  Given a family of hyperplanes $\mathcal{F}$ and a point $p \in \rr^d$, we define the depth of $p$ in $\mathcal{F}$ as the minimum number of members of $\mathcal{F}$ that a ray starting from $p$ can hit.  This was introduced in \cite{RoH99}, and it was conjectured that every finite family $\mathcal{F}$ of hyperplanes in general position in $\rr^d$ has a point $p$ at depth greater than or equal to ${|\mathcal{F}|}/{(d+1)}$.  Theorem \ref{theorem-karasevhyperplanes} implies that the answer is affirmative when $|\mathcal{F}|/(d+1)$ is a prime power.

\begin{problem}
Does Theorem \ref{theorem-karasevhyperplanes} hold if $r$ is not a prime power?
\end{problem}

Another family of variations appear if we have a family of convex sets which are large (i.e., they have large volume, large diameter, many lattice points, etc...) and we want to partition them so that the intersection of the convex hull of the parts is also a large convex set.  We call these \textit{quantitative} versions of Tverberg's theorem.  Take for example the following Tverberg-type theorem for the diameter \cite{Sob16diam}.

\begin{theorem}
Let $r,d$ be positive integers and $\varepsilon > 0$ a real number.  Then, there is a number $M = M(r,d,\varepsilon,\operatorname{diam})$ such that the following holds.  Given a family $X$ of $M$ intervals of length $1$ in $\rr^d$, there is a partition of $X$ into $r$ subfamilies $X_1, \ldots, X_r$ such that the diameter of $\cap_{j=1}^r \conv (\cup X_j)$ is at least $1-\varepsilon$.  Moreover, $M$ is linear in $r$.
\end{theorem}

The loss of diameter $\varepsilon$ is necessary for this result.  An equivalent statement can be proved for other functions \cite{RoS17quant}, such as the volume.  It is unclear if the loss $\varepsilon$ is still necessary in that setting.

\begin{problem}
Given $r,d$, determine if there is a number $M(r,d,\operatorname{vol})$ such that the following holds.  For any family $\mathcal{C}$ of $M$ convex sets in $\rr^d$, each of volume at least one, there is a partition of $\mathcal{C}$ into $r$ subfamilies $\mathcal{C}_1, \ldots, \mathcal{C}_r$ such that
\[
\operatorname{vol} \left(\bigcap_{j=1}^r \conv\left(\bigcup \mathcal{C}_j\right)\right) \ge 1.
\]
\end{problem}

Other interpretations of quantitative Tverberg appear in \cite{discrete-tverberg17, de17continuous}.  In those results, we are given a family of $n$ ``large'' convex sets $K_1, \ldots, K_n$, and we seek a transversal $y_1 \in K_1, \ldots, y_n \in K_n$ that admits a (usual) Tverberg partition but where the intersection of the convex hulls of the parts is also ``large''.

We can also significantly change the convexity in the conclusion of Tverberg's theorem.  As mentioned in the introduction, if a subset $X$ of $\rr^d$ with $|X|=(r-1)(d+1)+1$ is in sufficiently general position, then for a partition of $X$ into $r$ sets $X=X_1\cup \ldots \cup X_r$ $1\le |X_j|\le d+1$ for every $j$, then $\bigcap_{j=1}^r \aff X_j$ is a single point.

Tverberg's theorem simply says that we can always find a partition such that, if $\{p\} = \bigcap_{j=1}^r \aff X_j$, the coefficients of the affine combination of $X_j$ that give $p$ are non-negative.  It turns out that sometimes we can prescribe some of those coefficients to be negative \cite{barany2016tverberg}.

\begin{theorem}
Let $X$ be a subset of $\rr^d$ of $(r-1)(d+1)+1$ points in sufficiently strong general position.  Let $M \subset X$ be a set of points such that $\conv (M) \cap \conv (X \setminus M) = \emptyset$.  Then, there is a partition of $X$ into $r$ parts $X_1, \ldots, X_r$ such that in the affine combinations that witness   $\bigcap_{j=1}^r \aff X_j \neq \emptyset$,
either
\begin{itemize}
	\item all the coefficients for $M$ are negative and all the coefficients for $X\setminus M$ are positive, or
	\item all the coefficients for $M$ are positive and all the coefficients for $X \setminus M$ are negative.
\end{itemize}
\end{theorem}

\begin{corollary}
Assume that, under the conditions of the previous theorem, $|M| < r$. Then in the affine combinations that witness $\bigcap_{j=1}^r \aff X_j \neq \emptyset$, all the coefficients for $M$ are negative and all the coefficients for $X\setminus M$ are positive.
\end{corollary}

However, if we require $k$ coefficients be negative, but we do not prescribe which $k$ points will carry the negative coefficients, the values of $k$ for which this is possible is an open problem.

\begin{problem}
Find all triples of integers $d,r,k$ for which the following holds. Given a subset $X$ of $\rr^d$ of $(r-1)(d+1)+1$ point in sufficiently general position, there is a partition of $X$ into $r$ parts $X_1, \ldots, X_r$ such that among the affine combinations that witness  $\bigcap_{j=1}^r \aff X_j \neq \emptyset$, exactly $k$ coefficients are negative.
\end{problem}

\section{Acknowledgments}

This work was partly supported by the National Science Foundation under Grant No. DMS-1440140 while the first author was in residence at the Mathematical Sciences Research Institute in Berkeley, California, during the Fall 2017 semester. The first author was also supported by Hungarian National Research, Development and Innovation Office Grants no K111827 and K116769. The authors would like to thank Fr\'ed\'eric Meunier, Uli Wagner, G\"unter M. Ziegler, and an anonymous referee for their careful revision and helpful comments.
\newcommand{\etalchar}[1]{$^{#1}$}
\providecommand{\bysame}{\leavevmode\hbox to3em{\hrulefill}\thinspace}
\providecommand{\MR}{\relax\ifhmode\unskip\space\fi MR }
\providecommand{\MRhref}[2]{%
  \href{http://www.ams.org/mathscinet-getitem?mr=#1}{#2}
}
\providecommand{\href}[2]{#2}

\smallskip

\noindent
Imre B\'ar\'any \\
\textsc{
Alfr\'ed R\'enyi Institute of Mathematics,\\
Hungarian Academy of Sciences\\
H-1364 Budapest Pf.~127  Hungary\\
and\\
Department of Mathematics\\
University College London\\
Gower Street, London, WC1E 6BT, UK}\\
\textit{E-mail address: }\texttt{barany.imre@renyi.mta.hu}
\medskip

\noindent
\noindent Pablo Sober\'on \\
\textsc{
Mathematics Department \\
Northeastern University \\
Boston, MA 02445, USA}\\
\textit{E-mail address: }\texttt{pablo.soberon@ciencias.unam.mx}

\end{document}